\documentclass[12pt]{amsart}
\usepackage{latexsym}
\usepackage{accents}
\usepackage{amssymb}
\usepackage{amsmath, amscd}
\usepackage{a4}
\usepackage{setspace}
\newtheorem*{theorem*}{Theorem}
\newtheorem{theorem}{Theorem}[section]
\newtheorem{lemma}[theorem]{Lemma}
\newtheorem{proposition}[theorem]{Proposition}
\newtheorem{corollary}[theorem]{Corollary}
\newtheorem{definition}[theorem]{Definition}

\newtheorem{notation}[theorem]{Notation}

\newcommand\rmd{{\mathrm{d}}}

\newcommand\ch{{\mathcal {H}}}

\newcommand\ca{{\mathcal {A}}}

\newcommand\cj{{\mathcal {J}}}

\newcommand\bfb{{\mathbf {B}}}
\newcommand\bfk{{\mathbf {K}}}

\newcommand\balpha{{\boldsymbol {\alpha}}}
\newcommand\bbeta{{\boldsymbol {\beta}}}
\newcommand\bomega{{\boldsymbol {\omega}}}
\newcommand\bgamma {{\boldsymbol {\gamma}}}

\newcommand\blambda {{\boldsymbol {\lambda}}}
\newcommand\bmu {{\boldsymbol {\mu}}}

\newcommand\rms{{\mathrm {s}}}

\newcommand\rmh{{\mathrm {h}}}
\newcommand\rma{{\mathrm {a}}}

\newcommand\fri{{\mathfrak {i}}}

\newcommand\tens{{\otimes_{\tau}}}

\newcommand{\bb}[1]{\mathbb{#1}}

\newcommand{\cl}[1]{\mathcal{#1}}

\begin{document}

\title{  $S$-numbers  of elementary operators on $C^*$-algebras}

\author{M. Anoussis, V. Felouzis and I. G. Todorov}
\address{Department of Mathematics, University of the Aegean,
832 00 Karlovasi -- Samos, Greece}
\email {mano@aegean.gr}

\address{Department of Mathematics, University of the Aegean,
832 00 Karlovasi -- Samos, Greece}
\email {felouzis@aegean.gr}
\address{Department of Pure Mathematics, Queen's University Belfast,
Belfast BT7 1NN, United Kingdom}
\email {i.todorov@qub.ac.uk}
\date{2 October 2008}
\begin{abstract}
We study  the $s$-numbers  of elementary operators acting on $C^*$-algebras.
The main results are the following:
If $\tau$ is any tensor norm
and $A,B\in\bfb(\cl H)$ are such that the sequences 
$s(A),s(B)$ of their singular numbers belong to a
stable Calkin space $\fri$ then the sequence of approximation numbers
of $A\otimes_{\tau} B$ belongs to $\fri$. 
If $\mathcal{A}$ is a $C^{*}$-algebra, $\mathfrak{i}$ is a stable
Calkin space, ${\rm s}$ is an $s$-number function, and $a_i, b_i \in
\mathcal{A},$ $ i=1,2,\dots,m$ are such that $s(\pi(a_i)),
s(\pi(b_i)) \in \mathfrak{i}$, $i=1,2,\dots,m$ for some faithful
representation $\pi$ of $\cl A$ then ${\rm s}\left(\sum_{i=1}^{m}
M_{a_i,b_i}\right)\in \mathfrak{i}$. The converse implication holds
if and only if  the ideal of compact elements of $\cl A$
has finite spectrum. 
We also prove a quantitative  version of  a result of Ylinen. 
\end{abstract}
\maketitle


\section*{introduction}

Let $\cl A$ be a $C^{*}$-algebra. If $a,b\in\cl A$ 
we denote by  $M_{a,b}$
the operator on $\cl A$ given by $M_{a,b}(x) = axb$. An
operator $\Phi: \cl A\rightarrow \cl A$ is called {\it elementary}
if $\Phi = \sum_{i=1}^m M_{a_i,b_i}$ for some $a_i, b_i \in \cl
A$, $i=1,\dots,m$. 

Let $\cl H$ be a separable Hilbert space and $\bfb(\cl H)$
the $C^{*}$-algebra  of
all bounded linear operators on  $\cl H$.
A theorem of Fong and  Sourour \cite{fs} asserts that
an elementary operator $\Phi$ on  $\bfb(\cl H)$  is compact
if and only if there exists a representation $\sum_{i=1}^m
M_{A_i,B_i}$ of $\Phi$ such that the {\it symbols} $A_i, B_i$, $i =
1,\dots,m$ of $\Phi$ are compact operators. An element $a$ of a $C^*$-algebra $\cl A$ is called {\it compact} if the operator
$M_{a,a}$ is compact. Ylinen \cite{y} showed that $a\in \cl A$
is a compact element if and only if there exists a faithful
*-representation $\pi$ of $\cl A$ such that the operator $\pi(a)$
is compact.

The result of Fong and   Sourour was extended by Mathieu \cite{m} 
who showed that if $\cl A$ is a prime $C^*$-algebra, 
then an elementary
operator $\Phi$ on $\cl A$ is compact if and only if there exist
compact elements $a_i,b_i\in \cl A$, $i = 1,\dots,m$,
such that  $ \Phi=\sum_{i=1}^m
M_{a_i,b_i}$. Recently Timoney \cite{ti} extended this result 
to general $C^*$-algebras.

In this paper we investigate quantitative aspects of the above
results. It is well-known that a bounded
operator on a Banach space is compact
if and only if its  Kolmogorov numbers  form a null sequence. 
In our approach we use the more general
notion of the $s$-function introduced by  Pietsch
and the theory of ideals of $\mathbf{B}(\cl
H)$ developed by von Neumann, Schatten,  Calkin and others. A
detailed study of these notions is presented in the monographs
\cite{pi},\cite{carl}, \cite{gk} and  \cite{si}.

In Section \ref{sni} of the paper we recall the definitions of Calkin spaces
and the basic properties of $s$-functions. 

In Section
\ref{ideal_sets} we study  stable Calkin spaces.
An analogous property for ideals of $\mathbf{B}(\cl H)$,  called ``tensor product
closure property'', was 
considered by  Weiss \cite{wei}.
We give a necessary and sufficient condition for the stability of
a singly generated Calkin space. We also provide a sufficient condition for the stability of  a Lorentz sequence space.

If $a, b \in \cl A$ and  $\cl C$ is a $C^*$-subalgebra
of $\cl A$ such that $M_{a,b}(\cl C)\subseteq \cl C $ 
 we  denote by $M_{a,b}^{\cl C}$ the operator $\cl C \rightarrow \cl C$ defined by $ M_{a,b}^{\cl C}(x)=axb$. 
In Section \ref{restrict}  we prove inequalities relating $s$-number functions of the operators $ M_{a,b}$ and $ M_{a,b}^{\cl C}$.

In Section \ref{elbh} we study elementary operators acting on
$\mathbf{B}(\cl H)$. Some of our results can be presented in a more
general setting. Namely, we show that if $\tau$ is any tensor norm
and $A,B\in\bfb(\cl H)$ are such that $s(A),s(B)$ belong to a
stable Calkin space $\fri$ then the sequence of approximation numbers
of $A\otimes_{\tau} B$ belongs to $\fri$. 
A result of this type
for $\fri = \ell_{p,q}$ was proved by K\"{o}nig in \cite{ko} who
used it to prove results concerning tensor stability of $s$-number
ideals in Banach spaces. 
We also show that if $\Phi$ is an elementary
operator on $\bfb (\cl H)$, $\mathfrak{i}$ is a stable Calkin space
and $\mathrm{s}$ is an $s$-function then $\mathrm{s}(\Phi)\in
\mathfrak{i}$ if and only if there exist $A_i,
B_i\in \bfb(\cl H)$, $i = 1,\dots,m$, such that  $\Phi=\sum_{i=1}^m
M_{A_i, B_i}$ and $s(A_i), s(B_i)\in \mathfrak{i}$. It is well known that all $s$-functions coincide for operators acting
on Hilbert spaces.
 It follows from our result that if $\Phi$ is an elementary operator on
$\mathbf{B}(\cl H)$,  $\fri$ is a stable Calkin space and $\mathrm{s}$,
 $\mathrm{s}'$  are  $s$-number functions, then the 
sequence $\mathrm{s}(\Phi)$ belongs to $\fri$ if and only if 
 the sequence $\mathrm{s'}(\Phi)$ belongs to
$\mathfrak{i}$.

In Section \ref{s_calg} we study elementary operators acting on 
$C^*$-algebras. We show  that if
$\mathcal{A}$ is a $C^{*}$-algebra, $\mathfrak{i}$ is a stable
Calkin space, ${\rm s}$ is an $s$-number function, and $a_i, b_i \in
\mathcal{A},$ $ i=1,\dots,m$, are such that $s(\pi(a_i)),
s(\pi(b_i)) \in \mathfrak{i}$, $i=1,\dots,m$ for some faithful
representation $\pi$ of $\cl A$ then ${\rm s}(\sum_{i=1}^{m}
M_{a_i,b_i})\in \mathfrak{i}$. The converse implication holds
if and only if  the ideal of compact elements of $\cl A$
has finite spectrum. 
Finally,  we  prove that  if
 $a\in\cl A$ and 
$\mathrm{d}(M_{a, a}) \in \mathfrak{i}$ for some Calkin space
$\mathfrak{i}$ then $s(\rho(a))^2\in \mathfrak{i}$, where $\rho$ is the reduced atomic representation of 
$\cl A$. 
This result may be viewed as a quantitative  
version of the aformentioned result of Ylinen.


\section{calkin spaces and $s$-functions}\label{sni}

In this section we  recall some notions and results
concerning the ideal structure of the algebra of all bounded
linear operators acting on a separable Hilbert space.
We also  recall the definition of an  $s$-function.

We will denote by $\bfb$ the class of all bounded linear operators
between Banach spaces. If $\cl X$ and $\cl Y$ are Banach spaces,
we will denote by $\bfb(\cl X, \cl Y)$ the set of all bounded
linear operators from $\cl X$ into $\cl Y$.
If $\cl X=\cl Y$ we set  $\bfb(\cl
X)=\bfb(\cl X, \cl X)$. Ideals of $\bfb(\cl X)$ or, more
generally, of a normed algebra $\ca$, will be proper, two-sided
and not necessarily norm closed. By $\mathbf{K}(\cl X)$ (resp.
$\mathbf {F}(\cl X)$) we  denote the ideal of all compact
(resp. finite rank) operators on $\cl X$. By $\|T\|$ we denote the operator
norm of a bounded linear operator $T$. We denote by $\ell_\infty$
the space of all bounded complex sequences, by $c_0$ the space of
all sequences in $\ell_{\infty}$ converging to $0$ and by $c_{00}$
the space of all sequences in $c_0$ that are eventually zero. The
space of all $p$-summable complex sequences is denoted by
$\ell_p$.
 For a subspace $ \jmath$ of $\ell_{\infty}$,
we let $\jmath^+$ be the subset of $\jmath$ consisting of all
sequences with non-negative terms. We denote by $c_0^\star$ the
subset of $c_0^+$ consisting of all non-negative decreasing
sequences.

If $\balpha = (\alpha_n)_{n=1}^\infty$ and $\bbeta =
(\beta_n)_{n=1}^\infty$ are sequences of real numbers, we write
$\balpha\leq\bbeta$ if $\alpha_n\leq \beta_n$ for each
$n\in\bb{N}$.
 For every $\balpha=(\alpha_n)_{n=1}^\infty\in c_0$ we define
$\balpha^\star=(\alpha_n^\star)_{n=1}^\infty\in c_0^\star$ to be
the sequence given by
$$
\begin{array}{lllll}
\alpha^\star_1=\max\{|\alpha_n|: n\in \mathbb{N}\},\\[2ex]
\alpha_1^\star+\cdots+\alpha_n^\star= \max\left\{\sum_{i\in
I}|\alpha_i|: I\subseteq \mathbb{N}, |I|=n\right\}.
\end{array}
$$
The sequence $\balpha^\star$ is the rearrangement of the sequence
$(|\alpha_n|)_{n=1}^\infty$ in decreasing order  including
multiplicities.

A\textit{ Calkin space} \cite{si} is a subspace $\mathfrak{i}$
of 
$c_0$ which has the following property:

\begin{enumerate}
 \item [] If $\balpha\in \fri$ and $\bbeta\in c_0$ then 
$\bbeta^\star\leq \balpha^*$ implies that $\bbeta\in \fri$.
\end{enumerate}

Let $\ch$ be a separable Hilbert space. If $e,f\in \ch$ we denote
by $f^*\otimes e$  the rank one operator on $\cl H$ 
given by $f^*\otimes
e(x) = (x,f)e$, $x\in \ch$. Given an operator $T\in \bfk(\ch)$
there exist orthonormal sequences $(f_n)_{n=1}^\infty\subseteq \cl
H$ and $(e_n)_{n=1}^\infty\subseteq \cl H$ and a unique sequence
$(s_n(T))_{n=1}^\infty\in c_0^\star$ such that
$$T = \sum_{n=1}^\infty s_n(T) f^*_n\otimes e_n$$ where the series converges in  norm. Such a decomposition of $T$ is called a \textit{Schmidt
expansion}. The elements of the sequence
$s(T)=(s_n(T))_{n=1}^\infty$ are called \textit{singular numbers} of $T$.

For every ideal $\mathcal{I}\subseteq \bf  B(\cl H)$ we set
$$\mathfrak{i}(\mathcal{I})=\{
\balpha\in c_0: \; \text{there exists}\;\; T\in\mathcal{I}
\;\text{such that}  \; \balpha^\star = s(T)\};$$ conversely, for
every Calkin space $\mathfrak{i}$ we set
$$\mathcal{I}(\mathfrak{i})=\{ T\in \bfb(\cl H): s(T)\in \mathfrak{i}\}.$$

The following classical result  describes the ideal structure of
$\mathbf{B}(\cl H)$ in terms of Calkin spaces (for a proof
of the formulation given here
 see \cite[Theorem 2.5]{si}).

\medskip

\noindent {\bf Theorem }{\bf (Calkin \cite{cal})} {\it Let $\cl H$ be a separable infinite dimensional
Hilbert space. The mapping $\mathcal{I} \mapsto
\mathfrak{i}(\mathcal{I})$ is a bijection from the set of 
all ideals of $\bfb(\cl H)$ onto the set  of all Calkin spaces 
 with inverse $\mathfrak{i}\mapsto
\mathcal{I}(\mathfrak{i})$.}

\medskip

We now recall  Pietsch's definition  of $s$-functions.
 A map $\mathrm{s}$ which assigns to every
operator $T\in \bfb$ a sequence of non-negative real numbers
$\rms(T)=(\rms_1(T), \rms_2(T), \dots)$ is called an
\textit{$s$-function} if the following are satisfied:

\begin{enumerate}
 \item  $\|T\|=\rms_1(T)\geq \rms_2(T)\geq \dots$, for $T\in \bfb$.
\item $\rms_n(S+T)\leq \rms_n(S)+\|T\|$, for $S, T\in \bfb(\cl X,
\cl Y)$.
\item $\rms_n(RST)\leq \|R\|\; \|T\|\rms_n(S)$, for $T\in \bfb(\cl
X, \cl Y)$, $S\in \bfb(\cl Y, \cl Z)$, $R\in \bfb(\cl Z, \cl W)$.
\item If $\mathrm{rank}(T)<n$ then $\rms_n(T)=0$.
\item $\rms_n(I_n)$=1, where $I_n$ is the identity operator on
$\ell_2^n$. (Here  $\ell_2^n$ is the $n$-dimensional 
complex Hilbert space).
\end{enumerate}

\noindent An $s$-function $\rms$  is said to be \textit{additive}
 if $\rms_{m+n-1}(S+T)\leq
\rms_m(S)+\rms_n(T)$  for all $m, n$ and  $S,T\in \bfb(\cl X, \cl Y)$.

We give below the definition of some $s$-functions which will be used
in the sequel. Let $\cl X$ and $\cl Y$ be Banach spaces and
$T\in\bfb(\cl X,\cl Y)$.

\medskip

\noindent (a)\;\;\; The sequence
$\mathrm{a}(T)=(\mathrm{a}_n(T))_{n=1}^\infty$ of
\textit{approximation numbers} of $T$ is given by
$$\mathrm{a}_n(T)=\inf\left\{\|T-A\| : A\in \bfb
(\cl X,\cl Y),
\ \mathrm{rank}(A) < n\right\}.$$

\medskip

\noindent(b)\;\;\;The sequence
$\mathrm{d}(T)=(\mathrm{d}_n(T))_{n=1}^\infty$ of
\textit{Kolmogorov numbers} of $T$ is given by
$$ \mathrm{d}_n(T)=\inf_{V}\;\;\sup_{\|x\|\leq 1}
\;\;\inf_{y\in V} \|Tx- y\|,$$ where the infimum is taken over all
subspaces $V$ of $\cl Y$ with $\dim V < n$.

\medskip

\noindent(c)\;\;\;The sequence $\mathrm{h}(T)$ of \textit{Hilbert
numbers} of $T$ is given by 
$$\mathrm{h}_n(T)=\sup s_n(ATB)$$
 where
the supremum is taken over all contractions $A\in \bfb(\cl
Y,\ch)$, $B\in \bfb(\cl K,\cl X)$ and Hilbert spaces $\cl H$ and
$\cl K$.

\medskip
The approximation, the Kolmogorov  and the Hilbert $s$-functions are additive \cite{pi}. 
Moreover,
for every $s$-function $\rms$, every
operator $T\in \bfb$ and every $n$ we have
$\rmh_n(T) \leq \rms_n(T) \leq \rma_n(T)$ \cite{pi}.

A well-known result of Pietsch (\cite[Theorem 11.3.4 ]{pi})
implies that if $\mathrm{s}$ is an $s$-function, $\cl H$ is a
separable Hilbert space and $T\in \bfk(\ch)$ then
$\mathrm{s}_n(T)$ is equal to the $n^{th}$-singular number
$s_n(T)$ of $T$.


\section{stable calkin spaces}\label{ideal_sets}

In this section we present some results concerning stable Calkin spaces.
 We characterize the stable principal  
Calkin spaces and show that certain Lorentz sequence
spaces are stable Calkin spaces.

If $\balpha=(\alpha_n)_{n=1}^\infty,
\bbeta=(\beta_n)_{n=1}^\infty\in c_0$, we define the sequence
$\balpha\otimes \bbeta = (\gamma_n)_{n=1}^\infty\in c_0^\star$ by
$$
\begin{array}{lllll}
\gamma_1=\max\{|\alpha_i\beta_j|: (i, j)\in \mathbb{N} \times \mathbb{N}\}\\[2ex]
\gamma_1 + \cdots+\gamma_n = \max\left\{\sum_{(i, j)\in
I}|\alpha_i\beta_j|: I\subseteq \mathbb{N} \times \mathbb{N},
|I|=n\right\}.
\end{array}
$$
The sequence $\balpha\otimes \bbeta$ is the rearrangement of the
double sequence $(|\alpha_n \beta_m|)_{n, m=1}^\infty$ in decreasing
order including multiplicities.

\begin{definition}\label{sta}
Let $\mathfrak{i}$ and $\mathfrak{j}$ be Calkin spaces. We let
$\mathfrak{i}\otimes \mathfrak{j}$ be the smallest Calkin space
containing the sequences $\balpha\otimes\bbeta$, where
$\balpha\in\mathfrak{i}$ and $\bbeta\in\mathfrak{j}$. A Calkin space
$\mathfrak{i}$ is said to be \textrm{stable} if
$\mathfrak{i}\otimes \mathfrak{i}=\mathfrak{i}$.
\end{definition}

Let $\cl H$ be a separable infinite dimensional Hilbert space. 
Weiss \cite{wei} defined the  \textit{tensor product closure property}
for ideals of $\bf B(\cl H)$. An ideal $\mathcal{I}$ of $\bf B(\cl H)$
has  this property if
$S\otimes T\in \mathcal{I}$ whenever $S,T\in \mathcal{I}$. Here
the Hilbert space $\cl H\otimes \cl H$ is identified with $\cl H$
in a natural way. It is easy to see that an ideal
$\mathcal{I}\subseteq \bf B(\cl H)$ has the tensor product closure property  if and only if
$\fri(\cl I)$ is a stable Calkin space.
We will need the following lemma.

\begin{lemma}\label{mon}
If $\balpha,\balpha',\bbeta,\bbeta'\in c_0^+$ are such that
$\balpha\leq\balpha'$ and $\bbeta\leq\bbeta'$ then
$\balpha\otimes\bbeta\leq \balpha'\otimes\bbeta'$.
\end{lemma}
\begin{proof}
Clearly $\alpha_m\beta_n\leq \alpha_m'\beta_n'$
for every $m$, $n$. So to prove the lemma it  suffices
to prove that 
if $\balpha\leq \bbeta$ then 
$\balpha^\star\leq \bbeta^\star$.
Consider an injection
$\pi: \mathbb{N}\rightarrow \mathbb{N}$
such that $\alpha_n^\star=\alpha_{\pi(n)}$.
Clearly $\alpha_1^\star\leq \beta_1^\star$.
Suppose that  
$\alpha_i^\star\leq \beta_i^\star$ for $i=1, \dots, n-1$.
If $\alpha_n^\star> \beta_n^\star$ then $\alpha_i^\star > \beta_n^\star$ 
for every $i=1, \dots, n$  and hence 
$\beta_{\pi(1)}, \dots, \beta_{\pi(n)}$ are 
strictly greater than 
$\beta_n^\star$ so the number 
of all $i\in \mathbb{N}$ such that 
$\beta_i>\beta_n^\star$ is greater than $n-1$,
a contradiction. 
\end{proof}

\begin{notation}
\begin{enumerate}
 \item If $\balpha_n=(\alpha^n_i)_{i=1}^{k_i}$ are finite
sequences, we set
$$(\balpha_n)_{n=1}^\infty = (\alpha^1_1, \dots, \alpha_{k_1}^1, \alpha_1^2, \dots,
\alpha^2_{k_2}, \alpha^3_1, \dots ,\alpha^3_{k_3},\ldots).$$
\item\ If $\omega\in\bb{C}$ and $r\in \mathbb{N}$ we set
$(\omega)_r =(\underbrace{\omega,\dots,\omega}_{r \; times})$.
\item If $r\in\bb{N}$ we let $\mathbf{r}=(\underbrace{1, 1, \dots,
1}_{r\;\; times}, 0,0,\dots).$
\item  \ If $\balpha, \bbeta$ are sequences of real numbers we write
$\balpha \lesssim \bbeta$ if there exists a constant $C>0$ such
that $\balpha\leq C\bbeta$.
\end{enumerate}

\end{notation}

\noindent Observe that if $(M_n)_{n=0}^\infty$, 
$(N_n)_{n=0}^\infty$ are sequences of non negative integers
then

\begin{equation}\label{eq2}
((\omega^n)_{M_n})_{n=0}^\infty\;\otimes\; ((\omega^n)_{N_n})_{n=0}^\infty=
((\omega^n)_{\sum_{k+l=n}M_k N_l})_{n=0}^{\infty}.
\end{equation}

\begin{lemma}\label{re}
Let $\mathfrak{i}$ be a Calkin space. If $\balpha\in \mathfrak{i}$
and $\bbeta\in c_{00}$ then $\balpha\otimes \bbeta\in \mathfrak{i}$.
\end{lemma}\begin{proof}
Let $\balpha^{\star} = (\alpha_n^{\star})_{n=1}^{\infty}$.
Clearly, if $r\in \mathbb{N}$ then
$$\mathbf{r}\otimes \balpha
= ((\alpha_1^\star)_r, (\alpha_2^\star)_r, \dots,
(\alpha_n^\star)_r, \dots).$$ It  follows  from the
definition of a Calkin space that $\mathbf{r}\otimes \balpha \in
\mathfrak{i}$. Since $\bbeta^\star\in c^+_{00}$, there exists $r\in
\mathbb{N}$ such that $\bbeta^\star \lesssim \mathbf{r}$. By Lemma
\ref{mon}, $\bbeta\otimes\balpha = \bbeta^\star \otimes \balpha
\lesssim \mathbf{r}\otimes \balpha.$ Hence, $\bbeta \otimes
\balpha\in \mathfrak{i}$.
\end{proof}

The following notation will be used in the sequel.

\begin{notation}\label{not}
Let $\balpha = (\alpha_m)_{m=1}^\infty\in c_0^{\star}$  and
$\omega\in (0,1)$. 
For every $n=0,1,\dots$, set
\begin{align*}
&\cl K_n^{(\omega)}(\balpha)=\{ m: \omega^{n+1}<\alpha_m\leq \omega^n\}, \ \
K_n^{(\omega)}(\balpha)=|\cl K_n^{(\omega)}(\balpha)|,\\
&\widetilde{K}_n^{(\omega)}(\balpha)=\sum_{i=0}^{n}
K_i^{(\omega)}(\balpha), \ \
M_n^{(\omega)}(\balpha)=\sum_{i+j=n}K_i^{(\omega)}(\balpha)
K_j^{(\omega)}(\balpha),\\
&\widetilde{M_n}^{(\omega)}(\balpha)=\sum_{i=0}^n
M_i^{(\omega)}(\balpha), \\ & K_{-1}^{(\omega)}(\balpha)=\widetilde{K}_{-1}^{(\omega)}(\balpha)=M_{-1}^{(\omega)}(\balpha)=\widetilde{M}_{-1}^{(\omega)}(\balpha)=0.
\end{align*}

\end{notation}
\begin{lemma} \label{lem2}
Let $\omega\in (0,1)$, $\balpha=(\alpha_m)_{m=1}^\infty\in c_0^{\star}\;,
\bbeta=(\beta_m)_{m=1}^\infty\in c_0^{\star}$. 
Assume that $\alpha_1, \beta_1\leq 1$. Then $\balpha\lesssim
\bbeta$ if and only if there exists a positive integer $r$ such that
for every $n \in \mathbb{N}\cup \{0\}$,
$$\widetilde{K}_{n}^{(\omega
)}(\balpha)\leq \widetilde{K}_{n+r}^{(\omega
)}(\bbeta).$$
\end{lemma}
\begin{proof}
Set $\widetilde{K}_n=\widetilde{K}_{n}^{(\omega )}(\balpha)$ and
$\widetilde{L}_n=\widetilde{K}_{n}^{(\omega )}(\bbeta)$. Suppose
that $\balpha\lesssim \bbeta$ and let $C > 0$ be such that
$\alpha_m\leq C \beta_m$, for every $m\in\bb{N}$. Let $r\in\bb{N}$ be such
that $\omega^r C\leq 1$. Then $\beta_{\widetilde{K}_n}\geq
C^{-1}\alpha_{\widetilde{K}_n}
> C^{-1} \omega^{n+1}\geq \omega^{n+1+r}$. Thus, $\widetilde{K}_n\leq \widetilde{L}_{n+r} $.

Conversely, suppose that there exists $r \in \mathbb{N}$ such that $\widetilde{K}_n\leq \widetilde{L}_{n+r}$, for every 
$n\in\bb{N}\cup\{0\}$. Fix
$m\in\bb{N}$ and let $n$ and $k$ be such that $m =
\widetilde{K}_{n-1}+ k$ and $1\leq k\leq  K_n$. Since $\widetilde
K_{n-1}< m\leq \widetilde{K}_{n}\leq \widetilde{L}_{n+r} $ we have

$$\alpha_m\leq \omega^{n} =
\omega^{-r-1} \omega^{n+r+1}\leq \omega^{-r-1}
\beta_{\widetilde{L}_{n+r}} \leq \omega^{-r-1} \beta_m. $$ Thus,
$\balpha\lesssim \bbeta$.
\end{proof}

If $\balpha\in c_0$ we let $\langle\balpha\rangle$ denote the smallest Calkin space
 containing $\balpha$. A Calkin space of the form $\langle\balpha\rangle$
is called \textit{principal}.
The proof of the following lemma is staightforward.
\begin{lemma}\label{lempr}
If $\balpha\in c_0$ then $$\langle\balpha\rangle=\{\bbeta\in c_0: \; \text{there exists} \;r\in \mathbb{N} \; \text{such that} \;\bbeta^\star\lesssim \mathbf{r}\otimes \balpha\}.$$
\end{lemma}

\begin{theorem}\label{propstable}
Let $\balpha=(\alpha_n)_{n=1}^\infty\in c_0^{\star}$ 
with $\alpha_1\leq 1$ and $\omega\in (0,1)$. The following
are equivalent:
\begin{enumerate}
 \item The principal Calkin space $\langle\balpha\rangle$ is stable.
\item$\balpha\otimes \balpha\in \langle\balpha\rangle$.
\item There exists $r\in \bb{N}$ and $C > 0$ such that
$\widetilde{M}_n^{(\omega)}(\balpha)\leq C
\widetilde{K}_{n+r}^{(\omega)}(\balpha)$, for every $n\in\bb{N}\cup\{0\}$.
\end{enumerate}

\end{theorem}
\begin{proof}
(1)$\Rightarrow$(2) is trivial.

(2)$\Rightarrow$(1) Let $\bbeta, \bgamma\in \langle\balpha\rangle$. By
Lemma \ref{lempr}, there exist positive integers $m, n$
such that $\bbeta\lesssim \bf m \otimes \balpha$ and
$\bgamma\lesssim \bf n \otimes \balpha$. By Lemma \ref{mon} we
have that $\bbeta\otimes\bgamma\lesssim (\bf m \otimes \balpha)
\otimes(\bf n \otimes \balpha)= (\bf m \bf n)\otimes
(\balpha\otimes\balpha)$. Since $\balpha\otimes \balpha\in
\langle\balpha\rangle$, using Lemma \ref{lempr}  again we conclude
that $ \bbeta\otimes\bgamma\in \langle\balpha\rangle$ and so $\langle\balpha\rangle$ is
stable.

(1)$\Leftrightarrow$(3) Set $K_n=K_n^{(\omega)}(\balpha)$ and  $\widetilde \balpha=(
(\omega^n)_{K_n} )_{n=0}^\infty$; clearly,
$\langle\balpha\rangle=\langle\widetilde\balpha\rangle$. By the previous paragraph,
$\langle\balpha\rangle$ is stable if and only if $\widetilde\balpha\otimes
\widetilde\balpha \in \langle\widetilde \balpha\rangle$. By  
Lemma \ref{lempr}, $\widetilde\balpha\otimes \widetilde\balpha \in
\langle\widetilde\balpha\rangle$ if and only if there exists a positive
integer $m$ such that $\widetilde\balpha\otimes \widetilde\balpha
\lesssim \bf m \otimes \widetilde\balpha$. Since $\widetilde
K_n^{(\omega)}(\widetilde\balpha\otimes \widetilde\balpha)=
\widetilde M_n^{(\omega)}(\widetilde\balpha)$ (see equation  (\ref{eq2}))
and $\widetilde K_{n}^{(\omega)}({\bf m} \otimes
\widetilde\balpha)= m \widetilde
K_{n}^{(\omega)}(\widetilde\balpha)$ the conclusion follows from Lemma \ref{lem2}.
\end{proof}

\begin{corollary}\label{remc} Let $\balpha=(\alpha_n)_{n=1}^\infty\in c_0^{\star}$ 
with $\alpha_1\leq 1$ 
and $\omega\in(0,1)$. Suppose that $C>0$ is
a constant such that
\begin{equation}\label{eqstt}
K^{(\omega)}_{n+j}(\balpha)\geq  C \left(\sum_{i=0}^n K^{(\omega)}_i(\balpha)\right)^2
\end{equation}
for all $j \in \mathbb{N}$ and $n\in\bb{N}\cup\{0\}$. Then $\langle\balpha\rangle$ is a stable Calkin space.
\end{corollary}
\begin{proof}
Set $K_n=K^{(\omega)}_n(\balpha)$, $\widetilde{K}_n =
\widetilde{K}_n^{(\omega)}(\balpha)$ and $\widetilde{M}_n =
\widetilde{M}_n^{(\omega)}(\balpha)$. Let $r$ be a positive integer
such that $rC\geq 1$. Since
$(\sum_{i=0}^nK_i)^2\geq\widetilde{M}_n$, it follows that 
$$\widetilde{K}_{n+r}\geq K_{n+1}+\dots+K_{n+r}\geq rC
\widetilde{M_n}\geq \widetilde{M_n}, \ \ n\in\bb{N}$$ and hence
condition (3) of Theorem \ref{propstable} holds.
\end{proof}

We next give some examples of stable and non-stable principal
Calkin spaces.

\bigskip

\noindent {\bf Examples}
\begin{enumerate}
 \item It follows from assertion (3) of Theorem
\ref{propstable} that for every $\omega\in (0,1)$, the principal
Calkin space $\left\langle(\omega^n)_{n=0}^\infty\right\rangle$ is not stable. This example
was first given in \cite{wei}.
\item\noindent Let $\lambda>0$ and $\balpha
=(n^{-\lambda})_{n=1}^\infty$. Then the principal Calkin space
$\langle\balpha\rangle$ is not stable. To show this, let $ \mu=\lambda^{-1}$
and $\omega=e^{-1}$. Let $K_n, M_n, \widetilde {K}_n,
\widetilde{M}_n$ be the positive integers associated with the
sequence $(n^{-\lambda})_{n=1}^\infty$ and $\omega$ (Notation
\ref{not}). Since $$\frac 12 \left [e^{(j+1)\mu}-e^{j\mu}\right
]\leq K_j\leq \left [e^{(j+1)\mu}-e^{j\mu}\right ], $$ there exist
constants $C_1, C_2
>0$ such that for every $j$
$$C_2 e^{j\mu}\leq K_j\leq C_1 e^{j\mu}. $$ 
It follows that
$M_n=\sum_{i+j=n}K_i K_j\geq (n+1)C_2^2e^{n\mu}$. Let $r$ be a positive integer. Then
$$\widetilde{K}_{n+r}\leq C_1\sum_{i=0}^{n+r} (e^\mu)^i=\frac{C_1(e^{(n+r+1)\mu}-1)}{e^{\mu}-1} $$
and
$$ \widetilde{M}_n\geq C_2^2\int_{0}^n (t+1)e^{\mu t}dt\geq
C_3n e^{n\mu}$$ for some $C_3> 0$. Thus,
$$\lim_{n\rightarrow+\infty}\frac{\widetilde{M}_n}{\widetilde{K}_{n+r}}=+\infty$$
for each $r\in\bb{N}$. By Theorem \ref{propstable}, $\langle\balpha\rangle$ is
not stable.

It follows from the characterization of the symmetrically 
 normable principal ideals due to Allen and
Shen \cite{all} that the principal ideal 
$\langle T\rangle$ of $\mathbf{B}(\cl H)$
generated by an operator $T$ with $s(T)=
 (n^{-\lambda})_{n=1}^\infty$,
$\lambda\in (0, 1)$ is symmetrically normed.
However, as we have shown, the principal Calkin space   $\left\langle(n^{-\lambda})_{n=1}^\infty\right\rangle$, for
$\lambda\in (0, 1)$,    is 
not  stable.
\item Let $\balpha=\left(\frac 1 {\log_2
m}\right)_{m=2}^\infty$. Then the Calkin space $\langle\balpha\rangle$ is  stable. To see
this,  consider the positive integers $K_n, M_n, \widetilde
{K}_n, \widetilde{M}_n$ associated with the sequence  $\left(\frac
1 {\log_2 m}\right)_{m=2}^\infty$ and $\omega = \frac{1}{2}$
(Notation \ref{not}). We have that
$K_n=2^{2^{n+1}}-2^{2^n}$.
Since
$$(K_0+\dots +K_n)^2=\left(2^{2^{n+1}}-2\right )^2$$
it follows from  Corollary \ref{remc} that
$\langle\balpha\rangle$ is stable.
\end{enumerate}

\medskip

In the sequel we  examine the stability of a class of Calkin spaces, namely,  the Lorentz
sequence spaces. We recall their definition \cite{lt}. Let $1\leq
p < \infty$ and let $\mathbf{w}=(w_n)_{n=1}^\infty$ be a
decreasing sequence of positive numbers such that $w_1=1$,
$\lim_{n\rightarrow \infty}w_n=0$ and $\sum_{n=1}^\infty
w_n=\infty$. We shall call  such a $\mathbf{w}$ a {\it
weight sequence}. The linear space $\ell_{\mathbf{w}, p}$ of all
sequences $\balpha=(\alpha_n)_{n=1}^\infty$ of complex numbers such that 

$$\|\balpha\|_{\mathbf{w}, p}\stackrel{def}{=}
\sup_{\pi}\left\{ \left(\sum_{n=1}^\infty w_n\;|\alpha_{\pi(n)}|^p
\right )^\frac 1p\right\}<\infty,$$ where $\pi$ ranges over all
the permutations of $\mathbb{N}$, is a Banach space under the
previously defined norm, called a {\it Lorentz sequence space}. 

If $\balpha\in\ell_{\mathbf{w}, p}$ then we easily see that
$$\|\balpha\|_{\mathbf{w}, p} =
\left( \sum_{n=1}^\infty w_n(\alpha_n^\star)^p\right)^\frac
1p.$$ If $w_n=n^{\frac pq-1}$ with $0< p<q$ we obtain the classical
$\ell_{q, p}$ spaces of Lorentz.

\begin{theorem}\label{lorenz}
Let $\mathbf{w}=(w_n)_{n=1}^\infty$ be a weight sequence such that
there exists a constant $C>0$ with $w_{mn}\leq C w_m w_n$ for
every $m, n\in \mathbb{N}$. Then for every $p\geq 1$ and $\balpha,
\bbeta \in \ell_{\mathbf{w}, p}$ we have that
$$\| \balpha \otimes \bbeta \|_{ \mathbf{w}, p} \leq
C^{1/p} \;\|\balpha\|_{\mathbf{w},p}\; \|\bbeta\|_{\mathbf{w},
p}.$$ In particular, $\ell_{\mathbf{w}, p}$ is a stable Calkin space.
\end{theorem}
\begin{proof}
We may assume  that $\balpha=(\alpha_n)_{n=1}^\infty$ and 
$\bbeta=(\beta_n)_{n=1}^\infty$ are  positive  decreasing sequences
with $\alpha_1, \beta_1\leq 1$. Fix $\omega$ with $0<\omega< 1$. For
every $n \in \mathbb{N}\cup\{0\}$ let $K_n=|\cl K_n^{(\omega)}(\balpha)|$, $L_n
= |\cl K_n^{(\omega)}(\bbeta)|$, $\widetilde M_n=\sum_{0\leq i+j\leq
n}K _i L_j$ and $K_{-1} = L_{-1} = \widetilde M_{-1}=0$.

Let
$$\tilde \balpha= \left((\omega^{n})_{K_n}\right)_{n=0}^\infty,\;\;\;
\tilde \bbeta= \left((\omega^{n})_{L_n}\right)_{n=0}^\infty.$$
Then
$$\tilde \balpha\otimes \tilde \bbeta=\left( (\omega^n)_{\widetilde M_{n}-\widetilde M_{n-1}}  \right)_{n=0}^\infty.$$
For every $n,i,k, l  \in \mathbb{N}\cup\{0\}$ such that  $0\leq i\leq n$,
$1\leq k\leq K_i$ and $1\leq l \leq  L_{n-i}$ we set
\begin{equation}\label{eqqq1}
\phi_n(i,k, l)=\widetilde M_{n-1}+\sum_{j=0}^{i-1}K_j L_{n-j} + kl.
\end{equation}
Also, for every $i$, $1\leq k\leq K_i$ and $1\leq l\leq L_i$ we set
$$\psi(i, k)=\sum_{j=0}^{i-1} K_j +k, \;\;\;\; \psi'(i, l)=\sum_{j=0}^{i-1} L_j +l.$$
We observe that for every positive integer $r$,  $(\tilde
\balpha\otimes \tilde \bbeta)_r = \omega^n$ if and only if
$r=\widetilde M_{n-1}+s$ with $1\leq s\leq \sum_{i+j=n} K_i L_j$ and
therefore $(\tilde \balpha\otimes \tilde \bbeta)_r=\omega^n$ if and
only if there exist $n,i,k, l  \in \mathbb{N}\cup\{0\}$ such that  $0\leq i\leq n$, $ 1\leq k\leq K_i$,
$1\leq l \leq L_{n-i}$  and $r=\phi_n(i, k, l)$. So,
\begin{equation}\label{8}
\|\tilde\balpha\otimes \tilde\bbeta\|^p_{\mathbf{w}, p} =
\sum_{n=0}^\infty \left( \sum_{i=0}^{n}\sum_{k=1}^{K_i}\sum_{l=1}^{
L_{n-i}} w_{\phi_n(i, k, l)} \right) \omega^{np}.
\end{equation}
Also,
$\tilde\alpha_r=\omega^i$ if and only if $r=\sum_{j=-1}^{i-1} K_j+k$ for some $k$ with
$1\leq k\leq  K_i$ and $\tilde\beta_r=\omega^{i'}$ if and only if $r=\sum_{j=-1}^{i'-1}L_j+l$ for some $l$ with
$1\leq l\leq L_{i'}$. So,

$$\|\tilde\balpha\|_{\mathbf{w}, p}^p=\sum_{n=0}^\infty \left(\sum_{k=1}^{K_n}w_{\psi(n, k)}
\right)\omega^{np}, \ \ \|\tilde\bbeta\|_{\mathbf{w},
p}^p=\sum_{n=0}^\infty \left(\sum_{l=1}^{L_n}w_{\psi'(n, l)}
\right)\omega^{np}$$ and
\begin{equation}\label{9}
\left \|\tilde\balpha \right \|_{\mathbf{w}, p}^p\;
\|\tilde\bbeta\|_{\mathbf{w}, p}^p\; =\sum_{n=0}^\infty \left(
\sum_{  i=0}^{n}\sum_{k=1}^{K_i}\sum_{ l=1}^{ L_{n-i}   } w_{\psi(i,
k)} w_{\psi'(n-i, l)}\right)\omega^{np}.
\end{equation}

\noindent But
$$
\begin{array}{lll}
\psi(i, k)\; \psi'(n-i, l)=\left(\sum_{j=0}^{i-1} K_j +k\right)\left(\sum_{j=0}^{n-i-1} L_j +l \right)=\\[3ex]
\sum_{j=0}^{i-1} \sum_{j'=0}^{n-i-1}K_j  L_{j'} +k \sum_{j=0}^{n-i-1} L_j +l \sum_{j=0}^{i-1} K_j+kl\leq\\[3ex]
\sum_{j=0}^{i-1} \sum_{j'=0}^{n-i-1}K_j  L_{j'} +K_i \sum_{j=0}^{n-i-1} L_j +L_{n-i} \sum_{j=0}^{i-1} K_j+kl\\[3ex]
\leq \widetilde M_{n-1}+kl \leq \phi_n(i, k, l).
\end{array}
$$
By the monotonicity of the weight sequence $\mathbf{w}$ we have
\begin{equation}\label{10}
w_{\phi_n(i, k, l)}\leq w_{\psi(i, k)\; \psi'(n-i, l)} \leq C
w_{\psi(i, k)}\; w_{\psi'(n-i, l)}.
\end{equation}
Finally, by (\ref{8}), (\ref{9}) and (\ref{10}), 
\begin{eqnarray*}
\|\balpha\otimes \bbeta\|_{\mathbf{w}, p} & \leq & \|\tilde
\balpha\otimes \tilde \bbeta\|_{\mathbf{w}, p} \leq C^{1/p}
\|\tilde\balpha\|_{\mathbf{w},p}\; \|\tilde\bbeta\|_{\mathbf{w}, p}\\
& = & C^{1/p} \dfrac 1 {\omega^2}\; \|\omega
\tilde\balpha\|_{\mathbf{w},p}\; \|\omega
\tilde\bbeta\|_{\mathbf{w}, p}\leq C^{1/p}\dfrac 1 {\omega^2}
\;\|\balpha\|_{\mathbf{w}, p}\; \|\bbeta\|_{\mathbf{w}, p}.
\end{eqnarray*}
Letting $\omega\rightarrow 1$ we obtain
$$\| \balpha \otimes \bbeta \|_{p, \mathbf{w}} \leq  C^{1/p} \;\|\balpha\|_{\mathbf{w},p}\;
\|\bbeta\|_{\mathbf{w}, p}.$$
\end{proof}


\section{$s$-numbers of restrictions  }\label{restrict}

Let $\mathcal{A}$ be a $C^*$-algebra.
 If $a,b\in\cl A$ 
we denote by  $M_{a,b}$
the operator on $\cl A$ given by $M_{a,b}(x) = axb$. An
operator $\Phi: \cl A\rightarrow \cl A$ is called {\it elementary}
if $\Phi = \sum_{i=1}^m M_{a_i,b_i}$ for some $a_i, b_i \in \cl
A$, $i=1,\dots,m$. 

 If  $\cl C$ is a $C^*$-subalgebra
of $\cl A$ such that $M_{a,b}(\cl C)\subseteq \cl C $ 
 we will denote by $M_{a,b}^{\cl C}$ the operator $\cl C \rightarrow \cl C$ defined by $ M_{a,b}^{\cl C}(x)=axb$. 
In this section we prove inequalities concerning $s$-number functions of the operators $ M_{a,b}$ and $ M_{a,b}^{\cl C}$.

 It is well-known that every closed
two-sided ideal $\cj$ of $\ca$ is an $M$-ideal, that is, that there
exists a  projection $\eta: \ca^*\rightarrow \cj^\perp$,
where $\cj^\perp$ is the annihilator of $\cj$ in $\ca^*$, such
that for every $\varphi\in \ca^*$,
\begin{equation*}
\|\varphi\| = \|\eta(\varphi)\| + \|\varphi-\eta(\varphi)\|
\end{equation*}
(see e.g. \cite{dav}, Theorem 11.4). A functional $\varphi\in
\ca^*$ is called a \textit{Hahn-Banach extension} of $\phi\in
\cj^*$ if it is an extension of $\phi$ and $\|\varphi\|=\|\phi\|$.
 If $\cj$ is an $M$-ideal of $\ca$ then
every $\phi\in \cj^*$ has  a unique Hahn-Banach  extension in
$\ca^*$  denoted by $\tilde{\phi}$. Thus, if we
identify $\cj^*$ with the subspace  $\{\widetilde\phi: \phi\in \cj^*\}$ of $\cl A^{*}$
then
\begin{equation*}
\ca^*=\cj^*\oplus_{\ell_1}\cj^{\perp};
\end{equation*}
thus $\|\widetilde \phi +\psi\|=\|\phi\|+\|\psi\|$ for all
$\phi\in \cj^*$, $\psi\in \cl J^\perp$. Given $T\in \bf B(\cj)$
let $\widehat{T}: \ca^*\rightarrow \ca^*$ be given by
$$\widehat{T}(\tilde\phi + \psi)=\widetilde{T^*(\phi)},$$ where
$\phi\in \cl J^*$ and $\psi \in \cl J^{\perp}$. We identify $\cl
A$ with a subspace of $\cl A^{**}$ via the canonical embedding and
denote by $\widetilde T : \cl A \rightarrow \cl A^{**}$ the restriction of
$\widehat{T}^*$ to $ \ca$.

\begin{lemma}\label{nei}
\begin{enumerate}
 \item 
 Let $T\in \bfb(\cl J)$. Then the operator $\widetilde{T}$ extends
 $T$ and $\|\widetilde{T}\| = \|T\|$.
\item
 The map $T\rightarrow \widetilde{T}$ is linear.
\end{enumerate}
\end{lemma}
\begin{proof}
The second assertion is easily verified. We show (1).  Let $x\in\cl J$ and $f\in \ca^*$. Then $f = \widetilde{\phi} +
\psi $ with $\phi \in \cj^*$ and $\psi \in \cj^{\perp}$. We have
\begin{eqnarray*}
\widetilde{T}(x)(f) & = & \widehat{T}^*(x)(f) =
\widehat{T}(f)(x) = \widetilde{T^*(\phi)}(x) = T^*(\phi)(x)\\ & =
& \phi(Tx) = \widetilde{\phi}(Tx) = f(Tx) = T(x)(f).
\end{eqnarray*}
Hence, $\widetilde{T}$ is an extension of $T$ and so
$\|T\|\leq\|\widetilde{T}\|$. 

We show that
$\|\widetilde{T}\| \leq \|T\|$.
Let $x\in\cl A$ and $f\in\cl A^*$. Then   $f = \widetilde{\phi} +
\psi $ with $\phi \in \cj^*$ and $\psi \in \cj^{\perp}$. We have
\begin{eqnarray*}
|\widetilde{T}(x)(f)| & = & |\widehat{T}^*(x)(f)| =
|\widehat{T}(f)(x)| = |\widetilde{T^*(\phi)}(x)| \leq
\|\widetilde{T^*(\phi)}\|\|x\|\\ & = & \|T^*(\phi)\| \|x\|\leq
\|T^*\| \|\phi\|\|x\|= \|T^*\|\|\widetilde{\phi}\|\|x\|\leq
\|T^*\|\|f\|\|x\|.
\end{eqnarray*}
Hence $\|\widetilde{T}\| \leq \|T^*\| = \|T\|$ and the proof is
complete.

\end{proof}

Let $\cl X$ be a reflexive Banach space and $\iota : \cl
J^*\rightarrow\cl A^*$ be the  map defined by  $\iota(\phi)
= \tilde{\phi}$, $\phi\in \cl J^*$. Clearly, $\|\iota\| = 1$.

Let $T : \cl J\rightarrow\cl X$. 
Write $T^{\sharp} : \cl A\rightarrow\cl X$ for the restriction of
$(\iota\circ T^*)^*$ to $\cl A$.

\begin{lemma}\label{lex}
The operator $T^{\sharp}$ extends $T$ and $\|T^{\sharp}\| =
\|T\|$.
\end{lemma}
\begin{proof}
Let $a\in \cl J$ and $g\in \cl X^*$. We have that
$$\begin{array}{lll}
T^{\sharp}(a)(g)  = & (\iota\circ T^*)^*(a)(g) = a((\iota( T^*(g)) = \iota( T^*(g))(a) =\widetilde{T^*(g)}(a)\\ [1ex]& =T^*(g)(a) =  g(Ta) =
T(a)(g).
\end{array}$$
Hence  $T^{\sharp}$ extends $T$ and so $\|T\|\leq \|T^{\sharp}\|$. On the other hand,
$$\|T^{\sharp}\| \leq \|(\iota\circ T^*)^*\| = \|\iota\circ T^*\|
\leq \|\iota\|\|T^*\| = \|T\|.$$
\end{proof}

\begin{lemma}\label{doh}
Let $\cl A$ be a $C^*$-algebra, $\cl J\subseteq\cl A$ be a closed two
sided ideal and $\Phi : \cl A\rightarrow\cl A$ be a bounded
operator which leaves $\cl J$ invariant. Let $\Phi_0:\cl J\rightarrow \cl J $ be the operator defined by $\Phi_{0}(x)=\Phi(x)$.  Then $h_n(\Phi_0)\leq h_n(\Phi)$
for each $n\in\bb{N}$.
\end{lemma}
\begin{proof} Write $\iota_0 : \cl J\rightarrow\cl A$ for the inclusion
map. In the supremum below, $\cl H$ and $\cl K$ are arbitrary
 Hilbert spaces. Using Lemma \ref{lex} we have that

$\begin{array}{lll}
h_n(\Phi_0) & = & \sup\{s_n(A\Phi_0 B) : B\in {\bf B}(\cl H,\cl J),
A\in {\bf B}(\cl J,\cl K) \mbox{ contractions}\}\\[0.8ex]
& = & \sup\{s_n(A^{\sharp}\Phi (\iota_0 B)) : B\in {\bf B}(\cl H,\cl
J), A\in {\bf B}(\cl J,\cl K) \mbox{ contractions}\}\\[0.8ex]
& \leq & \sup\{s_n(A_1\Phi B_1) : B_1\in {\bf B}(\cl H,\cl A), A_1\in
{\bf B}(\cl A,\cl K) \mbox{ contractions}\}\\[0.8ex]
& = & h_n(\Phi).\end{array}$\\
\end{proof}

If $\cl X$ is a Banach space,  $c \in \cl X$ and  $\phi \in \cl X^*$ we denote by
$\phi\otimes c$ the operator on $\cl X$ defined by $\phi\otimes c(x)=\phi(x)c$. We denote by $\mathbf{F}_n(\cl X)$ the set of all operators $F$ on $\cl X$ of rank less than or equal to $n$. It is well-known that $\mathbf{F}_n(\cl X)=\left\{\sum_{i=1}^{n}\phi_{i}\otimes c_i:\ \phi_{i} \in \cl X^{*}, c_i \in \cl X, i=1,2,\dots, n \right\}.$

\begin{lemma}\label{l_emb}
Let $\ca$ be a $C^*$-algebra and $\cj$ be a closed two-sided ideal of
$\ca$. 
\begin{enumerate}
 \item 
 Assume that $a, b\in \cj$.  Then  \ $\widetilde{M_{a, b}^\cj}(x)=M_{a,b}(x)$ for every $x \in \cl A$.
\item
 Let $\phi_i\in \cj^*, c_i\in \cj$, $i = 1,\dots,n$. Let $F$ be the operator on $\cl J$ defined by  
$F=\sum_{i=1}^n \phi_i\otimes c_i$.  Then
$\widetilde{F}(x)=\left(\sum_{i=1}^n \widetilde\phi_i\otimes c_i\right)(x)$ for every $x \in \cl A$.
\end{enumerate}
\end{lemma}
\begin{proof}
(1) Let $S = M_{a,b}$, $T = M_{a,b}^\cj$, and $\phi\in\cl J^*$.
First note that $S^*(\widetilde{\phi})$ is an extension of
$T^*(\phi)$. Indeed, for every  $x\in\cl J$ we have that 
$$S^*(\widetilde{\phi})(x) = \widetilde{\phi}(Sx) = \phi(Tx) =
T^*(\phi)(x).$$We show that $S^*(\widetilde{\phi})$ is the 
Hahn-Banach extension of $T^*(\phi)$. To this end, let $x\in\cl
A$  and $\{u_{\lambda}\}_{\lambda\in\Lambda}\subseteq
\cl J$ be a contractive approximate unit for $\cj$.  Then for each
$x \in \ca$, \ $au_{\lambda}xb \rightarrow_{\lambda} axb$
in norm and hence
$\phi(au_{\lambda}xb)\rightarrow_{\lambda} \phi(axb)$.
We thus have that
\begin{eqnarray*}
|S^*(\widetilde{\phi})(x)| & = & |\widetilde{\phi}(Sx)| =
|\widetilde{\phi}(axb)| = |\phi(axb)| = \lim_{\lambda\in\Lambda}
|\phi(au_{\lambda}xb)|\\ & = &\lim_{\lambda\in\Lambda}
|T^*(\phi)(u_{\lambda}x)| \leq \|T^*(\phi)\|\|x\|.
\end{eqnarray*}
It follows that $\|S^*(\widetilde{\phi})\|\leq \|T^*(\phi)\|$.
Since $S^*(\widetilde{\phi})$ extends $T^*(\phi)$, we have that
$S^*(\widetilde{\phi})$ is the Hahn-Banach extension of $T^*(\phi)$,
that is, $S^*(\widetilde{\phi}) = \widetilde{T^*(\phi)}$.

Let $x\in\cl A$ and $f\in\cl A^*$. Then  $f = \widetilde{\phi} +
\psi $ with $\phi \in \cj^*$ and $\psi \in \cj^{\perp}$ and
\begin{eqnarray*}
\widetilde{T}(x)(f) & = & \widehat{T}^*(x)(f) =
\widetilde{T^*(\phi)}(x) = S^*(\widetilde{\phi})(x) =
\widetilde{\phi}(Sx)\\ & = & (\widetilde{\phi} + \psi)(Sx) =
S^*(f)(x) = S(x)(f).
\end{eqnarray*}

(2) By Lemma \ref{nei}(2) it suffices to show the statement in the case $F =
\phi_{1}\otimes c_{1}$, where $\phi_{1}\in \cj^*$ and $c_{1}\in \cj$. Let $\phi \in\cl
J^*$. We have  that $\widetilde{F^*(\phi)}(x) =
\widetilde{\phi}_{1}(x)\phi(c_{1})$ for every $ x\in\cl A$. Indeed, the functional
$x\rightarrow \widetilde{\phi}_{1}(x)\phi(c_{1})$  extends $F^*(\phi)$
and has norm equal to the norm of $F^*(\phi)$ since $\|\phi_{1}\| =
\|\widetilde{\phi}_{1}\|$. Let $x\in\cl A$ and $f\in\cl A^*$. We have $f = \widetilde{\phi}
+ \psi$, where $\phi \in \cl J^*$ and $\psi\in \cl J^{\perp}$.  Then
\begin{eqnarray*}
\widetilde{F}(x)(f) & = & \widehat{F}^*(x)(f) = \widehat{F}(f)(x)
=
\widehat{F}(\widetilde{\phi} + \psi)(x) = \widetilde{F^*(\phi)}(x)\\
& = & \widetilde{\phi}_{1}(x)\phi(c_1) =
\widetilde{\phi}_{1}(x)\widetilde{\phi}(c_{1}) = \widetilde{\phi}_{1}(x)f(c_{1}) =
(\widetilde{\phi}_{1}\otimes c_{1})(x)(f).
\end{eqnarray*}
\end{proof}

The following theorem is the main result of this section.

\begin{theorem}\label{lem55}
Let $\cl A$ be a $C^*$-algebra, $\cl J$ be a closed two-sided ideal
of $\cl A$ and $a, b\in \cl J$. Then for every $n\in \bb{N}$ we
have that
$$
\mathrm{h}_n\left(M_{a,b}^{\cl J}\right)\leq
\mathrm{h}_n(M_{a,b})\leq \mathrm{a}_n\left(M_{a,b}\right)\leq
\mathrm{a}_n(M_{a,b}^\cj).
$$
\end{theorem}
\begin{proof}
The first inequality follows from Lemma \ref{doh} while the second
one is trivial. In what follows the operators $\widetilde{F}$ for $F\in \mathbf{ F}_{n-1}(\cl J)$ and $\widetilde{M_{a,b}^{\cj}}$ are considered as operators from $\cl A$  to $\cl A$; this is possible by Lemma \ref{l_emb}. It follows from Lemmas \ref{nei} and
\ref{l_emb} that for every $n \in \mathbb{N}$ we have
$$
\begin{array}{llll}
\mathrm{a}_n\left(M_{a,b}\right)
&= \inf \left \{\left \|M_{a,b}-G\right \|:\;\; G\in \mathbf{ F}_{n-1}(\ca)\right \}\\[2ex]
&\leq \inf \left \{ \left \|M_{a,b}-\widetilde{F}\right \|:\;\;  F\in \mathbf{ F}_{n-1}(\cl J) \right \}\\[2ex]
&= \inf \left \{\left \|\widetilde{M_{a,b}^{\cj}}-\widetilde{F}\right \|:\;\;  F\in \mathbf{F}_{n-1}(\cl J) \right \}\\[2ex]
&=\inf \left \{\left \|M_{a,b}^\cj-F \right \|: \;\;  F\in
\mathbf{ F}_{n-1}(\cl J) \right \} = \mathrm{a}_n(M_{a,b}^\cj).
\end{array}$$
\end{proof}

We close the section with a lemma  which will be used in the proof of Theorem \ref{th_main_ver}.
\begin{lemma}\label{l_kap}
Let $\cl B\subseteq \bf  B(\cl H)$ be a $C^*$-algebra, $\cl A =
\overline{\cl B}^{\mathrm{wot}}$ and $A\in\cl A$. 
Assume that $A\in \cl B$. 
Then $\mathrm{d}(M_{A, A})\leq \mathrm{d}\left(M_{A, A}^{\cl B}\right)$.
\end{lemma}
\begin{proof} 
Set $\rmd\left(M^{\cl B}_{A, A}\right)=(d_n)_{n=1}^\infty$.
Let $\epsilon > 0$ and $\cl F\subseteq\cl B$ be a linear
space such that $\dim \cl F < n$ and
$$\inf_{F\in\cl F}\|AXA - F\| < d_n + \epsilon$$ for each contraction $X\in \cl B$.
It suffices to show that $\inf_{F\in\cl F}\|AYA - F\| \leq d_n +
\epsilon$ for each contraction $Y\in \cl A$.
 Suppose this is not the
case and let $Y\in \cl A$ be a contraction such that $\|AYA - F\| >
d_n + \epsilon$, for each $F\in \cl F$. By the Kaplansky Density
theorem, there exists a net $(X_{\nu})_{\nu}\subseteq\cl B$ of
contractions such that $X_{\nu}\rightarrow_{\nu} Y$ in the weak
operator topology. Let $F_{\nu}\in\cl F$ be such that $\|AX_{\nu} A
- F_{\nu}\| < d_n + \epsilon$. We have that $\|F_{\nu}\|\leq d_n +
\epsilon + 1$ for each $\nu$, and hence we may assume without loss
of generality that $F_{\nu}\rightarrow F_0$ in norm. We thus have
$AX_{\nu}A - F_{_\nu}\rightarrow AYA - F_0$ weakly. It follows that
$$\|AYA - F_0\|\leq\liminf \|AX_{\nu}A - F_{\nu}\| \leq d_n +
\epsilon,$$ a contradiction.
\end{proof}

\section{Elementary operators  on $\mathbf{B}(\cl
H)$}\label{elbh}

In this section we obtain estimates for the $s$-numbers of an elementary
operator acting on $\bfb(\cl H)$ in terms of the singular numbers of its
symbols. We formulate some of our results using tensor
products. Recall \cite{pi} that a {\it cross norm} $\tau$ is a
norm defined simultaneously on all algebraic tensor products $\cl
X\otimes \cl Y$ of Banach spaces $\cl X$ and $\cl Y$ such that
$\tau(x\otimes y)=\|x\| \|y\|$ for all $x\in \cl X$ and $y\in \cl
Y$. By $\cl X\tens \cl Y$ we denote the completion of the
algebraic tensor product with respect to $\tau$. A {\it tensor norm}
is a cross norm $\tau$ such that for every $A\in\bfb(\cl X,\cl Y)$
and $B \in\bfb(\cl X',\cl Y')$ the linear operator $A\otimes B :
\cl X\otimes\cl X'\rightarrow\cl Y\otimes\cl Y'$ given by
$A\otimes B (x\otimes x')=Ax\otimes Bx'$ is bounded with respect
to $\tau$ and the norm of its extension $A \tens B\in \bfb(\cl
X\tens \cl X',\cl Y\tens \cl Y')$ satisfies the inequality $\|A
\tens B\|\leq \|A\|\; \|B\|$.

In Theorem \ref{prop40} below we give an upper bound for the
approximation numbers of the operator $A\tens B$ in terms of the
sequence $ s(A)\otimes s(B) $. We will need the following lemma
which is due to K\"{o}nig \cite[Lemma 2]{ ko}.

\begin{lemma}\label{lemmastable}
Let $\tau$ be a tensor norm, $\cl X, \cl Y$ be Banach spaces,
$A\in \mathbf{B}(\ell_2, \cl X)$, $B\in \mathbf{B}(\ell_2, \cl
Y)$ and $(P_k)_{k=0}^n$, $(Q_k)_{k=0}^n$ be families of mutually
orthogonal projections acting on $\ell_2$. Then
$$\left\|\sum_{k=0}^n AP_k\tens B Q_k\right\|_{\ell_2\tens\ell_2\rightarrow \cl X\tens \cl Y}\leq
\max_{0\leq k\leq n}\{ \|A P_k\|\; \|B Q_k\|\}.$$
\end{lemma}

\begin{theorem}\label{prop40}
Let $\cl H$ be a  Hilbert space, $A,B\in \bfk(\cl H)$
and $\tau$ be a tensor norm. Then
\begin{equation}\label{eq}\rma(A\tens B)\leq  6.75  \;\; s(A)\otimes s(B) .
\end{equation}
Consequently, if $\mathfrak{i}$ and $\mathfrak{j}$ are  Calkin spaces,
$s(A)\in \mathfrak{i}$ and $s(B)\in \mathfrak{j}$ and $\mathrm{s}$
is any $s$-function then $\mathrm{s}\left(A\tens B\right)\in \mathfrak{i}\otimes \mathfrak{j}$.\end{theorem}
\begin{proof} If $\balpha = (\alpha_n)_{n=1}^\infty$ is a bounded sequence we write $D_\balpha\in\bfb(\ell_2)$
for the diagonal operator given by $D_\balpha((x_n)_{n=1}^\infty) = (\alpha_n x_n)_{n=1}^\infty$ for  $(x_n)_{n=1}^\infty \in \ell_2$.
It suffices to prove the theorem in the case
where $A$ and $B$ are diagonal operators in $\bfb(\ell_2)$. Indeed, suppose
that (\ref{eq}) holds in this case. By  polar decomposition, there
exist partial isometries $U_A, U_B: \cl H\rightarrow \ell_2$, $V_A,
V_B: \ell_2\rightarrow \cl H$ and diagonal operators $D_\balpha,
D_\bbeta: \ell_2\rightarrow \ell_2$, where $\balpha=s(A), \bbeta=s(B)$, such that $A= V_A D_\balpha U_A$ and $B=V_B D_\bbeta
U_B.$ Then
$$ A\tens B=(V_A\tens V_B)  (D_\balpha \tens D_\bbeta) (U_A\tens  U_B),\;\; $$
and hence
\begin{eqnarray*}
\rma\left( A\tens B\right) & \leq & \| V_A\tens V_B\| \;\rma
\left(D_\balpha\tens D_\bbeta\right)\; \|U_A\tens  U_B\|
\leq \rma \left(D_\balpha\tens D_\bbeta\right)\leq\\
& \leq & 6.75\; \balpha \otimes \bbeta= 6.75\; s(A)\otimes s(B).
\end{eqnarray*}
So let $A=D_\balpha: \ell_2\rightarrow \ell_2, \ B=D_\bbeta:
\ell_2\rightarrow \ell_2$, where $\balpha=(\alpha_n)_{n=1}^\infty,
\bbeta=(\beta_n)_{n=1}^\infty$ are non-negative decreasing sequences.
We may further assume  that $\alpha_1, \beta_1\leq 1$. Set
$\rma_n=\rma_n\left(A\tens B\right)$ and fix $\omega$ such that  $0<\omega<
1$. 

In what folows we use the notation introduced in \ref{not}. For every $n \in \mathbb{N}\cup\{0\}$ let
$$
\begin{array}{llllll}
K_n = K_n^{(\omega)}(\balpha), \ \ L_n = K_n^{(\omega)}(\bbeta), \ \ \widetilde M_n=\sum_{0\leq i+j\leq n}K _i L_j, \ \
\widetilde M_{-1}=0\\[2ex]
P_n = \sum_{i\in \mathcal{K}_n^{(\omega)}(\balpha)}e_i^*\otimes
e_i,\;\; Q_n=\sum_{i\in \cl K_n^{(\omega)}(\bbeta)}e_i^*\otimes
e_i,
\end{array}
$$
where $(e_n)_{n=0}^\infty$ is the standard basis of $\ell_2$.

Let
$A_n =  A P_n$, $B_n=B Q_n$ and $E_n=\sum_{0\leq k+l \leq n} A_k\tens  B_l.$
Clearly, $\|A_n\|\leq \omega^n$, $\|B_n\|\leq \omega^n$ and
$\mathrm{rank}  E_n\leq \widetilde M_n.$ Moreover,
$$ A=\sum_{n=0}^\infty A_n, \;\; B=\sum_{n=0}^\infty B_n, \;\;
A\tens B=\sum_{n, m=0}^\infty A_m\tens  B_n,$$ where the series
are absolutely convergent in the norm topology. Hence, $$ \rma_{\widetilde
M_n+1}\left(A\tens B\right)\leq \left\|A\tens B-E_n\right\|\leq
\sum_{N=n+1}^\infty\left\| \sum_{k+l=N} A_k\tens B_l\right\|.$$ By
Lemma \ref{lemmastable},
\begin{eqnarray*}
\left\| \sum_{k+l=N} A_k\tens B_l\right\| & = & \left\| \sum_{k=0}^N
A P_k\tens  B Q_{N-k}\right\|\leq \max _{0\leq k\leq N}\|A_k\|\;
\|B_{N-k}\|\\ & \leq & \max _{0\leq k\leq N}\omega^k
\omega^{N-k}=\omega^N
\end{eqnarray*}
and so
\begin{equation}\label{eq1}
\rma_{\widetilde M_n+1}\leq \sum_{N=n+1}^\infty
\omega^N=\dfrac{1}{1-\omega} \; \omega^{n+1}.\end{equation}

By the monotonicity of the approximation numbers, Lemma \ref{mon},
(\ref{eq2}) and (\ref{eq1}) we obtain
\begin{eqnarray*}
(\rma_n)_{n=1}^{\infty} & = & \left(
 (\rma_j)_{j=\widetilde M_n+1}^{\widetilde M_{n+1}}
\right)_{n=-1}^{\infty} \leq \left( (\rma_{\widetilde
M_{n}+1})_{\widetilde M_{n+1}-\widetilde M_n}
\right)_{n=-1}^{\infty}\\ & \leq & \dfrac{1}{1-\omega} \left(
 (\omega^{n+1})_{\widetilde M_{n+1}-\widetilde M_n}
\right)_{n=-1}^{\infty} = \dfrac{1}{1-\omega} \left(
 (\omega^{n})_{\sum_{i+j=n}K_i L_j}
\right)_{n=0}^{\infty}\\
& = & \dfrac{1}{1-\omega} ((\omega^n)_{K_n})_{n=0}^\infty \otimes
((\omega^n)_{L_n})_{n=0}^\infty\\ & = &
\dfrac{1}{\omega^2(1-\omega)}((\omega^{n+1})_{K_n})_{n=0}^\infty
\otimes  ((\omega^{n+1})_{L_n})_{n=0}^\infty\\
& \leq &  \dfrac{1}{\omega^2(1-\omega)} s(A)\otimes s(B).
\end{eqnarray*}
The minimal value of  $\frac{1}{\omega^2(1-\omega)}$ for $\omega\in (0,1)$  is 
$6.75$, and so 
\begin{equation*}
\rma \left(A\tens B \right)\leq  6.75\;s(A)\otimes s(B).
\end{equation*}
 \end{proof}

Theorems \ref{prop40} and \ref{lorenz} yield the following corollary.

\begin{corollary}\label{cor11}
Let $\mathbf{w}=(w_n)_{n=1}^\infty$
be a weight sequence with $w_{mn}\leq C w_m w_n$ for all $m$,$n$ and $A, B\in
\mathbf{K}(\cl H)$ be operators with $s(A),  s(B)\in
\ell_{\mathbf{w}, p}$. Then
$$\|\rma(A\tens B)\|_{\mathbf{w}, p}\leq 6.75 \; C^{1/p}\;
\|s(A)\|_{\mathbf{w}, p}\;\|s(B)\|_{\mathbf{w}, p}.$$
\end{corollary}

Consider the weight sequence $\mathbf{w}=(w_n)_{n=1}^\infty$, where $w_n=\frac{(1+\ln n)^\gamma} {n^{\alpha}}$.  If $0 < \alpha\leq
1$ and $\gamma\geq 0$, then $w_{mn}\leq  w_m w_n$ for all $m$,$n$. Hence Corollary \ref{cor11}
extends results of H. K\"{o}nig (\cite[Proposition 3]{ko}) and F.
Cobos and L. M. Fern\`{a}ndez-Cabrera \cite{cob}.

\bigskip

For the rest of the paper, we will be concerned with elementary
operators. Let $A, B$ be compact operators in $\mathbf{B}(\cl H)$. We recall that $M_{A, B}$ is the operator $\bfb(\cl H)\rightarrow\bfb(\cl H)$ defined by $M_{A, B}(X)=AXB$
and $M_{A, B}^{\bfk(\cl H)}$ is the operator $\bfk(\cl H)\rightarrow\bfk(\cl H)$ defined by $M_{A, B}^{\bfk(\cl H)}(X)=AXB$.
Theorems
\ref{lem55} and \ref{prop40} imply the following corollary.

\begin{corollary}\label{prop41}
Let $A,B$ be compact operators in $\bfb(\cl H)$. Then $$\rma(M_{A,B})\leq\rma(M_{A,B}^{\bfk(\cl H)})\leq  6.75  \;\; s(A)\otimes
s(B).$$
\begin{proof}
For every $x\in \cl H$ we denote by $f_x$ the functional on $\cl H$ defined by $f_x(y)={\langle y, x \rangle}$.
The   \textit{conjugate} \textit{space} $\bar{\cl H}$ of $\cl H$
is defined to be the set  $\{ f_x : x\in \cl H\}$ 
 with vector space operations
$ f_x+ f_y=f_{x+y},\; \lambda f_x=f_{\bar\lambda x}$
and   inner product
 given by $\langle f_x, f_y\rangle=\overline{\langle x, y\rangle}$.
For every $A\in \bfb(\cl H)$ we 
denote by $\bar A\in \bfb(\bar{\cl H})$
the operator defined by 
$ \bar A(f_x)=f_{Ax}$.

Note that the map $A\mapsto \bar A$ is a surjective 
conjugate linear isometry and that 
$s(A)=s(\bar A)$, for every compact operator $A$.

 Let $\epsilon$ be the injective tensor 
norm. The mapping 
$F: \bar{\cl H}\otimes H\rightarrow \bfb(\cl H)$
given by 
$F\left(\sum_{i=1}^n f_{x_i}\otimes y_i\right )=
\sum_{i=1}^n x_i^*\otimes y_i$
is a linear isometry (\cite[Ch. IV, Theorem 2.5]{ta})
 of $\bar{\cl H}\otimes_\epsilon \cl H$ onto 
$\bfk(\cl H)$.

We define 
$\tilde F: \bfb(\bar {\cl H}\otimes \cl H)\rightarrow \bfb(\bfk(\cl H))$
by $\tilde F(T)=F\circ T \circ F^{-1}$.
Clearly $\tilde F$ is a surjective linear isometry 
and  $\tilde F(T)$ is given by $\tilde F(T)(x^*\otimes y)=F(T(f_x\otimes y))$  for $x, y\in \cl H$.

For every  $\bar A\in \bfb(\bar{\cl H}) $, where $A\in \bfb(\cl H)$,
and every $B\in \bfb(\cl H)$ we have that 
\begin{equation}\label{eqtensor}
 \tilde F (\bar A\otimes_\epsilon B)=
M_{B^*, A}^{\bfk(\cl H)}.
\end{equation}

Indeed, for  every $x, y\in \cl H$,
$$ 
\begin{array}{lll}
 \tilde F (\bar A\otimes_\epsilon B)(x^*\otimes y)=&
F(\bar A\otimes_\epsilon B)(f_x\otimes y)=F(\bar Af_x\otimes By)=F(f_{Ax}\otimes By)=\\[2ex]
&(Ax)^*\otimes By=B (x^*\otimes y) A^*=M_{B, A^*}^{\bfk(\cl H)}(x^*\otimes y).
\end{array}
$$
So if $A, B\in \bfb(\cl H)$ by (\ref{eqtensor}) and Theorems
\ref{lem55} and \ref{prop40} 
we have that 
$$
\begin{array}{lll}
\mathrm{a}\left( M_{A, B}\right)&\leq
\mathrm{a}\left( M_{A, B}^{\bfk(\cl H)}\right)
=\mathrm{a}\left(\tilde F \left( \bar B\otimes_\epsilon A^*\right)\right)=
\mathrm{a}\left( \bar B\otimes_\epsilon A^*\right)\\[2ex]
&\leq 6.75 s(\bar B)\otimes s(A^*)=6.75 s(A)\otimes s(B). 
\end{array}
$$
\end{proof}

\end{corollary}

\begin{proposition}\label{pr_hilb}
Let  $\ca$  be a $C^*$-subalgebra of $\bfb(\cl H)$ such
that $\mathbf{K}(\cl H)\subseteq \ca $. Let $A_i,
B_i\in \ca$, $i=1, \dots, m$, and $\Phi=\sum_{i=1}^m M_{A_i,B_i}$.  If the operators $A_i$ $($resp. $B_i$ $)$,
$i=1,\dots, m$, are linearly independent then there exists
$r\in\bb{N}$ and a constant $C > 0$ such that for every $n$ and for
every $i=1,\dots,m$,
$$s_{rn-r+1}(A_i) \leq C \;\mathrm{h}_n(\Phi) \;(\text{resp.}\;\;s_{rn-r+1}(B_i) \leq C \;\mathrm{h}_n(\Phi)).$$
In particular, if $\mathfrak{i}$ is a Calkin space  and
$\mathrm{h}(\Phi)\in\mathfrak{i}$ then $s(A_i)\in \mathfrak{i}$ $($resp.
$s(B_i)\in \mathfrak{i}$$)$  for every $i=1, \dots, m$.
\end{proposition}
\begin{proof}
We will only consider the case where the operators $B_i$ , $i=1,\dots, m$, are linearly
independent. The other case can be treated similarly.

By \cite[ Lemma 1]{fs}, there exist $r\in\bb{N}$ and $\xi_i, \eta_i
\in \mathcal{H}$, $i = 1,\dots,r$, such that
$$\sum_{j=1}^r\langle B_i\eta_j,\xi_j\rangle =
\begin{cases}
1\ &\text{if }i=1\\
0 &\text{if }i=2, \dots, m \
\end{cases}.$$
Let $\phi_j : \mathcal{H}\rightarrow \ca$, $j=1, \dots, r$, be the
operators given by $\phi_j(\xi)=\xi_j^*\otimes \xi$, $\psi_j:\cl A
\rightarrow \mathcal{H}$, $j=1, \dots, r$, be the operators given by
$\psi_j(B)=  B\eta_j$ and
$$S = \sum_{j=1}^r \psi_j \circ \Phi \circ \phi_j=\sum_{i=1}^{m} \sum_{j=1}^{r} \psi_j \circ M_{A_i, B_i} \circ \phi_j.$$
For $\xi\in \cl H$ we have
\begin{eqnarray*}
(\psi_j \circ M_{A_i,B_i} \circ \phi_j)(\xi) & = & \psi_j (A_i
\phi_j(\xi) B_i)
=\psi_j (A_i(\xi_j^*\otimes \xi)B_i)\\
& = & \psi_j((B_i^*\xi_j)^*\otimes A_i\xi)=\langle
\eta_j,B_i^*\xi_j\rangle A_i\xi = \langle B_i\eta_j,\xi_j\rangle
A_i\xi
\end{eqnarray*}
and hence
$$S = \sum_{i=1}^m \left(\sum_{j=1}^r\langle B_i\eta_j,\xi_j\rangle\right)  A_i = A_1.$$
By the additivity of the singular numbers, we have that
$$s_{rn-r+1}(A_1) = s_{rn-r+1}(S) \leq \sum_{j=1}^r s_n(\psi_j \circ \Phi \circ \phi_j), \ \ n\in\bb{N}.$$
Let $C = r\;\max_{j=1, \dots, r} \|\psi_j\|\|\phi_j\|$. Then $s_n(\psi_j
\circ \Phi \circ \phi_j) \leq \|\psi_j\|\|\phi_j\| h_n(\Phi)$ and so
$s_{nr-r+1}(A_1)\leq    C \mathrm{h}_n(\Phi)$,
$n\in\bb{N}.$ 

Finally, by the monotonicity of $s$-numbers, we have
that
\begin{eqnarray*}
s(A_1)=(s_n(A_1))_{n=1}^\infty & = &
((s_{nr-r+1+k}(A_1))_{k=0}^{r-1})_{n=1}^\infty\\ & \leq &
((s_{nr-r+1}(A_1))_r)_{n=1}^\infty \leq C\
((\mathrm{h}_n(\Phi))_r)_{n=1}^\infty.
\end{eqnarray*}
If $\mathfrak{i}$ is a Calkin space and $\mathrm{h}(\Phi)\in \mathfrak{i}$,
Lemma \ref{re} implies that $((h_n(\Phi))_r)_{n=1}^\infty\in
\mathfrak{i}$. It follows that $s(A_1)\in \mathfrak{i}$.
\end{proof}

The following theorem is the main result of this section.

\begin{theorem}\label{thst}
Let $\Phi$ be an elementary operator on $\bfb (\cl H)$ $($resp. on $\mathbf{K}(\cl H)$$)$, $\mathfrak{i}$ be
a stable Calkin space and $\mathrm{s}$ be an $s$-function. Then
$\mathrm{s}(\Phi)\in \mathfrak{i}$ if and only if there exist
$m\in\bb{N}$ and $A_i, B_i\in \bfb(\cl H)$, $i = 1,\dots,m$, such
that $\Phi=\sum_{i=1}^m M_{A_i, B_i}$  and $s(A_i), s(B_i)\in
\mathfrak{i}$ for $i=1, \dots, m$.
\end{theorem}
\begin{proof}
We prove the Theorem in the case where $\Phi$ is an elementary operator on $\bfb (\cl H)$. The proof in the case where $\Phi$ is an elementary operator on $\bfk (\cl H)$ is similar.

Suppose that $\mathrm{s}(\Phi)\in
\mathfrak{i}$. Let $\Phi = \sum_{i=1}^m M_{A_i, B_i}$ be a
representation of $\Phi$ where $m$ is  minimal.
Then $A_i$ (resp. $B_i$), $i=1,\dots,m$, are linearly independent. Since $\mathrm{h}(\Phi)\leq \mathrm{s}(\Phi)$ we
have that $\mathrm{h}(\Phi)\in \mathfrak{i}$. By Proposition
\ref{pr_hilb}, $s(A_i), s(B_i)\in \mathfrak{i}$ for every $i=1,
\dots,m$.

Conversely, suppose that $\Phi = \sum_{i=1}^m M_{A_i, B_i}$ where
$s(A_i), s(B_i)\in \mathfrak{i}$ for every $i=1,\dots,m$.
Since $\mathfrak{i}$ is stable, Corollary \ref{prop41} implies that
$\mathrm{a}(M_{A_i,B_i}) \in \mathfrak{i}$. By the additivity of the
approximation numbers, $\mathrm{a}(\Phi)\in \mathfrak{i}$ and so $\mathrm{s}(\Phi)\in \mathfrak{i}$.
\end{proof}

Theorem \ref{prop41} provides an upper bound for the the
approximation numbers of $M_{A,B}$ in terms of the sequence
$s(A)\otimes s(B)$.  In the following proposition we obtain a lower
bound for the  Hilbert numbers of $M_{A,B}$ in terms of the sequence
$s(A)\otimes s(B)$.  
For $1\leq p < \infty$ we denote by $(\cl S_p, \|\;\|_p)$  the Schatten $p$-class, that is, the space   of all operators $A \in \bf K(\cl H)$ such that $s(A) \in \ell^{p}$, where the norm is given by $\|A\|_p=\left(\sum_{n=1}^\infty |s(A)|^p \right)^{\frac{1}{p}}$. 
If $\balpha=(\alpha_n)_{n=1}^\infty$  and $\bbeta=(\beta_n)_{n=1}^\infty$ are sequences of complex numbers we denote by $\balpha \bbeta$ the sequence $(\alpha_n \beta_n)_{n=1}^\infty$

\begin{proposition}\label{lemh}
Let  $A, B$ be compact operators in  $\mathbf{B}(\cl H)$. The following hold:
\begin{enumerate}
\item If $\blambda$ and $\bmu$ are sequences of unit norm in $\ell_4^+$
then $\mathrm{h}(M_{A,B})\geq (\blambda s(A)) \otimes (\bmu s(B))$.\\
\item  If $\blambda$ is a sequence of unit norm in $\ell_2^+$ then
$\mathrm{h}(M_{A,B})\geq (\blambda s(A)) \otimes s(B)$ and
$\mathrm{h}(M_{A,B})\geq s(A) \otimes (\blambda s(B))$.
\end{enumerate}
In particular,
\begin{equation}\label{eqh3}
\mathrm{h}_n(M_{A,B})\geq \frac{(s(A) \otimes\ s(B))(n)}{\sqrt n}.
\end{equation}
\end{proposition}
\begin{proof}
(1)
Let $A,B\in \bfb(\cl H) $ be compact operators
of norm one and $A^*=U|A^*|$ and $B=V|B|$ be the polar
decompositions of $A^*$ and $B$, respectively.  Let  $s(A)=(\alpha_n)_{n=1}^\infty, \; s(B)=(\beta_n)_{n=1}^\infty$ and 
$$|A^*| = \sum_{i=1}^{\infty}\alpha_i e_i^*\otimes e_i \;\; \mbox{and} \;\; |B| =
\sum_{j=1}^{\infty}\beta_j f_j^*\otimes f_j$$ be Schmidt expansions
of $|A^*|$ and $|B|$, respectively. Let $\cl K$ be the closed subspace of $\mathcal{S}_2$ spanned by the family $\{f_i^*\otimes e_j, i, j\}$ and 
$F:\cl K \rightarrow \mathbf{B}(\cl H)$ be the map given 
by $F(X)=UXV^*$. Clearly $\|F\|\leq 1$.

 Consider sequences
$\blambda = (\lambda_i), \bmu = (\mu_j)\in \ell_4^+$ of unit
norm and let $D_{\blambda},D_{\bmu}\in\bfb(\cl H)$
be the operators given by
$$D_{\blambda} =\sum_{i=1}^\infty \lambda_i e_i^*\otimes e_i,
\;\;\;\;\;
 D_{\bmu} = \sum_{j=1}^\infty \mu_j f_j^*\otimes f_j.$$
Let  $G : \bfb(\cl H)\rightarrow \cl K$ be the operator given by
$G(Y) =D_\blambda Y D_\bmu$. Since
$$\|D_{\blambda}YD_{\bmu}\|_2\leq
\|D_{\blambda}\|_4\|D_{\bmu}\|_4\|Y\| \leq \|Y\|$$ the operator $G$ is well defined and $\|G\|\leq 1$. The family
$\{f_i^*\otimes e_j, i, j\}$ is an orthonormal basis of $\cl K$ and
$$(G \circ M_{A, B} \circ F) (f_i^*\otimes e_j) = \lambda_j \alpha_j \mu_i \beta_if_i^*\otimes e_j.$$
It follows that
$$\mathrm{h}_n(M_{A,B})\geq s_n(G \circ M_{A, B}\circ F ) =
(\blambda s(A) \otimes\bmu s(B))(n)$$ and (1) is proved.
The proof of (2) is similar. 

We show inequality (\ref{eqh3}).
Let $s(A) \otimes s(B)=
(\nu_n)_{n\in\mathbb{N}}$ and $\pi:\mathbb{N}\rightarrow  \mathbb{N}\times \mathbb{N}$ ,
$\pi(n)=(i_n,j_n)$ be a bijection such that  $\nu_n=\alpha_{i_n}
\beta_{j_n}$. We  set $\blambda = (\lambda_i)_{i=1}^\infty,\; \bmu =(\mu_j)_{j=1}^\infty$
where \\[1ex]
$$\lambda_i=\left\{
\begin{array}{llll}
\frac{1}{\sqrt[4]{n}} & \mathrm{if} \;\; i \in \{ i_1, \dots , i_n
\}
\\ [2ex] \; 0 & \mathrm{if}\;\; i\not \in \{i_1,\dots,i_n\}
\end{array}\right.,\;\;
\mu_j=\left\{
\begin{array}{ll}
\frac{1}{\sqrt[4]{n}} & \mathrm{if}\;\; i\in \{j_1,\dots,j_n\}\\
[2ex] \; 0 & \mathrm{if}\;\; i\not \in \{j_1,\dots,j_n\}
\end{array}\right. .
$$
We have that
$(\blambda s(A) \otimes \bmu s(B))(k)=\dfrac{1}{\sqrt{n}}\nu_k$ for
every $k = 1,\dots,n$ and so $\mathrm{h}_n(M_{A,B})\geq
\frac{1}{\sqrt{n}}\nu_n.$

\end{proof}

It follows from Theorem \ref{thst}  that if the $s$-numbers of the symbols of
an elementary operator $\Phi$ belong to a stable Calkin space
$\mathfrak{i}$ then the $s$-numbers of  $\Phi$ belong also to
$\mathfrak{i}$. In what follows we show that this is not true
 without the assumption that
$\mathfrak{i}$ be  stable.

\begin{proposition}\label{cor}
Let $\omega\in (0,1)$ and $\mathfrak{i}$ be the principal Calkin space
generated by the sequence $\bomega =(\omega^{n-1})_{n=1}^\infty$.
Then there exists $A\in \bf B(\cl H)$ such that $s(A)\in \fri$ and
$\mathrm{h}(M_{A,A})\not \in \mathfrak{i}$.
\end{proposition}
\begin{proof} Let $A\in \bf B(\cl H)$ be such that $s(A)=\bomega$.
We will show that $\mathrm{h}(M_{A,A})\not \in \fri$. By Proposition \ref{lemh} it suffices to show that the sequence
$\balpha=\left (\frac{1}{n} (\bomega\otimes \bomega)(n) \right)_{n=1}^\infty$
does not belong to $\mathfrak{i}$, or  (by Lemma \ref{lempr})
 that  for every $r\in \mathbb{N}$,
$\balpha\not \lesssim \bf r\otimes \bomega$.

Suppose that there exist $r_0\in
\bb{N}$ and $C > 0$ such that $\balpha\leq C {\bf
r_0}\otimes\bomega$. Let $\balpha = (\alpha_n)_{n=1}^{\infty}$ and
${\bf r_0 }\otimes\bomega = (\beta_n)_{n=1}^{\infty}$. Then for
every $m$ we have that
$$\beta_{r_0 m}=\omega^{m-1}\;\; \mathrm{and}\;\;  \alpha_{\frac{m(m+1)}{2}}=
\frac{2}{m(m+1)}\omega^{m-1}.$$ 
So, if $r$ is an even positive integer and  $n(r)=\frac{rr_0(rr_0+1)}
{2}$  we have that
$$\frac{2}{rr_0(rr_0+1)}\omega^{rr_0-1}=\alpha_{n(r)}\leq C\beta_{n(r)}=C\omega^{\frac{r(r_0r+1)}{2}-1} ,$$
which leads to  a
contradiction.
\end{proof}


\section{Elementary operators  on $C^*$-algebras}\label{s_calg}

Let $\cl A$ be a $C^*$-algebra. Recall that an element $a\in \cl A$ is called \textit{compact} if
the operator $M_{a,a}: \cl A\rightarrow \cl A$ is compact. We
denote by $\cl {K}(\cl A)$  the closed two-sided ideal of all
compact elements of $\cl A$. The \textit{spectrum} of $\cl A$ is the set of unitary equivalence classes of non-zero irreducible representations of $\cl A$.
We will need  two  lemmas  which follow from  
\cite[ \textsection 5.5]{mu}.

\begin{lemma}\label{rep}
Let $(\rho,\cl H)=\left(\bigoplus_{i \in I} \rho_i,\;\;
\bigoplus_{i \in I}\cl{H}_i\right)$ be the reduced atomic
representation of $\ca$  where  $\{(\rho_i, \cl{H}_i),
i\in I\}$ is a maximal family of unitarily inequivalent irreducible 
representations of $\ca$. Let $J=\{i \in I: \rho_i(\cl{K}(\cl A))\neq \{0\}\}$. 
Let $\sigma_{i}$ be the restriction of $\rho_i$  to $\cl K(\cl A)$.
Then the representation   $\sigma =\left( \bigoplus_{i \in J}\sigma _i,\;\;\bigoplus_{i \in J}\cl{H}_i\right)$ is the  reduced atomic representation of  $\cl{K}(\cl A)$.
\end{lemma}

\begin{lemma}\label{spe}
Let $\cl A$ be a $C^*$-algebra  such that  $\cl A=\cl {K}(\cl A)$ and   $\sigma=\left( \bigoplus_{i \in J}\sigma _i,\;\;\bigoplus_{i \in J}\cl{H}_i\right)$ be the  reduced atomic representation of $\cl A$. Then $\cl A$ has finite spectrum if and only if $J$ is finite. In this case,   $\sigma(\cl A)=\sum_{i\in J} \mathbf{K}(\cl{H}_{i})$.
\end{lemma}

\begin{theorem}\label{th_main}
Let $\mathcal{A}$ be a $C^{*}$-algebra, $\mathfrak{i}$ be a stable
Calkin space and $\mathrm{s}$ be an $s$-function.
 Let  $\Phi$ be a compact elementary
operator on $\ca$.
\begin{enumerate}
 \item Suppose that
\begin{equation}\label{eq_mini} \Phi=
\sum_{i=1}^{m} M_{a_i,b_i} \ \ \ \ 
a_i, b_i\in \cl A,\;i=1,\dots,m,
\end{equation}
and that $\pi$ is a faithful representation of $\mathcal{A}$ such
that $s(\pi(a_i)), s(\pi(b_i)) \in \mathfrak{i}$, $i=1,\dots,m$.
Then $\mathrm{s}(\Phi)\in \mathfrak{i}$.
\item Suppose that $\cl{K}(\cl A)$ has finite spectrum and that
$\mathrm{s}(\Phi)\in \mathfrak{i}$. Then there exist a
representation $\sum_{i=1}^{m} M_{a_i,b_i}, \ \ 
a_i, b_i\in \cl A,\;i=1,\dots,m$, of $\Phi$ and a faithful
representation $\pi$ of $\mathcal{A}$ such that $s(\pi(a_i)),
s(\pi(b_i)) \in \mathfrak{i}$, $i=1,\dots,m$.
\end{enumerate}
\end{theorem}
\begin{proof}
(1) Since $\mathrm{s}_n(\Phi)\leq \mathrm{a}_n(\Phi)$ for each $n$,
it suffices to show that $\mathrm{a}(\Phi)\in \mathfrak{i}$. By the additivity of 
the approximation numbers  we
have that $\mathrm{a}_{nm - m + 1}(\Phi)\leq \sum\nolimits_{i=1}^m
\mathrm{a}_n(M_{a_i,b_i}).$ If $\mathrm{a}(M_{a_i,b_i})\in
\mathfrak{i}$ for each $i = 1,\dots,m$, Lemma \ref{re} 
implies that   $\mathrm{a}(\Phi)\in \mathfrak{i}$. Thus, we may
assume that $\Phi = M_{a,b}$, where $a,b\in\ca$.

Let $\pi : \ca\rightarrow \bf B(\cl H)$ be a faithful
representation such that $s(\pi(a)), s(\pi(b))\in \mathfrak{i}$.
Set $A=\pi(a)$ and $B=\pi(b)$. We denote by $M_{A,B}$ the
corresponding elementary operator acting on $\pi(\cl A)$. Clearly,
$A$ and $B$ are compact operators and
$\mathrm{a}(\Phi)=\mathrm{a}(M_{A,B}).$ Let $\cj =
\pi(\ca)\cap\bfk(\cl H)$. By Theorem \ref{lem55},
$$
\mathrm{a}_n\left(M_{A,B}\right)\leq \mathrm{a}_n(M_{A,B}^\cj), \ \text{for every} 
\ \ n\in\bb{N}.
$$
The $C^*$-algebra $\cj$ is  equal  to a
$c_0$-direct sum $\oplus_{i\in I}\cj_i$, where $\cj_i =\mathbb{C}I_{m_i}\otimes   \bfk(\cl
H_i)$
where $m_i$ is a positive integer,  $I_{m_i}$ is the the identity operator 
on a Hilbert space of dimension $m_i$
and  $\cl H_i$ is  a Hilbert space  \cite[ Theorem 1.4.5 ] {arv}.

Let $\Theta: \cj \rightarrow \bf K(\cl H)$ be the canonical injection
and $\Delta :  \bf K(\cl H)\rightarrow \cj$ the operator given by
$\Delta(X) =\sum_{i\in I}P_iXP_i$, where $P_i=I_{m_i}\otimes Q_i$
and $Q_i$  is the orthogonal
projection from   $\bigoplus_{i\in I}\cl H_{i}$ onto $\cl H_i$. We have that $M^{\cj}_{A,
B}=\Delta \circ M^{\bf  K(\cl H)}_{A,B} \circ \Theta$ where $M^{\bf
K(\cl H)}_{A,B}$ is the corresponding elementary operator acting
on $\bfk(\cl H)$. Thus,

$$\mathrm{a}_n\left(M^{\cj}_{A,B}\right)\leq \|\Delta\|
\;\mathrm{a}_n\left(M^{\bfk(\cl H)}_{A,B}\right) \;\|\Theta\| \leq
\mathrm{a}_n\left(M^{\bfk(\cl H)}_{A,B}\right).$$
By Corollary \ref{prop41}, $\mathrm{a}(\Phi )\in \mathfrak{i}$.

\medskip

(2) We identify $\cl A$ with $\rho(\cl A)$ where $(\rho, \cl H)=\left(\bigoplus_{i \in I} \rho_i,\;\;
\bigoplus_{i \in I}\cl{H}_i\right)$ is the reduced atomic representation of $\cl A$.  By \cite[Theorem 3.1]{ti}, there exist $A_{0 j},B_{0 j}\in\cl
{K}(\cl A)$, $j = 1,\dots,m$, such that $\Phi = \sum_{j=1}^m
M_{A_{0 j},B_{0 j}}$. Since $\mathrm{h}_n(\Phi)\leq
\mathrm{s}_n(\Phi)$ for each $n$, we may assume that $\mathrm{s} =
\mathrm{h}$. 
Let $\Phi_0: \cl K(\cl A)\rightarrow \cl K(\cl A)$ be the operator defined by $\Phi_{0}(X)=\Phi(X)$. By Lemma
\ref{doh}, $\mathrm{h}(\Phi_0)\in \mathfrak{i}$. Consequently,
the  $C^*$-algebra $\cl K(\cl A)$ and the operator 
$\Phi_0$ satisfy our assumptions. Thus we may 
assume that $\cl A = \cl K(\cl A)$.

 By Lemmas \ref{rep} and \ref{spe},  $\cl K(\cl A)=\oplus_{i\in I_{0}}
\bfk(\cl H_i)$ where $I_{0}$ is a finite subset of $I$. Let  $i \in I_{0}$. 
Clearly, $\mathbf K(\cl {H}_i)$ is invariant by $\Phi$.  Let 
  $\Phi_i: \mathbf K(\cl {H}_i)\rightarrow \mathbf K(\cl {H}_i)$ be the operator defined by $\Phi_{i}(X)=\Phi(X)$. The operator  $\Phi_i$ is an elementary
operator on $ \mathbf K(\cl {H}_i)$. By Theorem \ref{lem55},
$\mathrm{h}(\Phi_i)\in \mathfrak{i}$. By Theorem
\ref{thst}, there exists a representation $
\sum_{j=1}^{m_i} M_{A_{i j},B_{i j}}$  of $\Phi_{i}$   where $A_{i j}, B_{i j}\in
\bfk(\cl H_i)$ and $s(A_{i j}), s(B_{i j})\in
\mathfrak{i}$. Considering $A_{i j}$ and $B_{i j}$ as operators on
$\cl H$ we obtain that $\Phi = \sum_{i = 1}^k
\sum_{j=1}^{m_i}M_{A_{i j},B_{i j}}$ is a representation with the
required properties.
\end{proof}

Part (2) of Theorem \ref{th_main} does not hold if 
we  do not assume that 
 $\cl{K}(\cl A)$ has finite spectrum. In fact,
we have the following:

\begin{theorem}\label{corfin}
Let $\ca$ be a $C^* $-algebra. The following are equivalent:
\begin{enumerate}
\item ${\cl K}(\cl A)$ has finite spectrum.

\item  Let $\rms$ be an $s$-function,  $\fri$ be a stable Calkin space and $\Phi$ be a compact elementary operator  on
 $\cl A$. Assume that 
$\mathrm{s}(\Phi)\in \fri$. Then there exist a representation
$\sum_{i=1}^n M_{a_i,b_i}$ of $\Phi$ and  a faithful representation
$\pi$ of $\cl A$ such that $s(\pi(a_i)), s(\pi(b_i))\in \fri$ for
every $i=1,\dots, n$.
\end{enumerate}
\end{theorem}
\begin{proof}
The implication (1) $\Longrightarrow$ (2) follows from  Theorem \ref{th_main}.
We prove that (2) implies (1).
Suppose that $\cl{K}(\cl A)\neq \{0\}$
and that $\cl{ K}(\cl A)$ does not have finite spectrum. We will show that for every $p>2$
there exists an elementary operator $\Phi$ on $\cl A$ such that

(a) \ $\mathrm{a}(\Phi)\in \ell_p$, and

(b) whenever $\pi$ is a faithful representation of $\cl A$ and
$\Phi = \sum_{i=1}^n M_{c_i,d_i}$, $c_i,d_i\in\cl{K}(\cl A)$,  there exists $i$, $1\leq i\leq n$, such that 
$s(\pi(c_i))\not\in\ell_p$ or $s(\pi(d_i))\not\in\ell_p$.

Let $\sigma$ be the reduced atomic representation of  $\cl{K}(\cl A)$. Then $$\sigma (\cl{K}(\cl A)) =\bigoplus_{j \in J} {\mathbf {K}}(\cl H_j).$$ It follows from Lemma \ref{spe} that $J$ is infinite. Choose an infinite
countable subfamily $\{\cl H_j\}_{j=1}^{\infty}$ of the family $J$.
 For each $j \in \mathbb{N}$, consider a unit vector $e_j\in  \cl H_j$. 

Let $r_j$ be the projection of $\cl K(\cl A)$ such that 
$\sigma( r_j)=e_j^*\otimes e_j$ and  $(\lambda_j)_{j=1}^{\infty}$ 
be a decreasing sequence of positive
real numbers belonging to $\ell_{2p}$ but not to $\ell_p$.
We set $c=\sum_{j=1}^\infty \lambda_j r_j$, 
$p_k = \sum_{j=1}^k  r_j$ and $\Phi = M_{c,c}\in \bf B(\cl A)$.
We will show that $\rma(\Phi)\in \ell_p$.

Let $\rho$ be  the reduced atomic representation of $\cl A$. 
  Let $c_{n}=\sum_{i=1}^{n} \lambda_i r_i$.   It follows from  Lemma \ref{rep} that 
$$M_{\rho(c_{n}),\rho(c_{n})}(\rho(a))=
\sum_{i=1}^{n}\sigma_{i}(r_{i})
\rho_{i}(a)\sigma_{i}(r_{i})$$ and hence  
the operator $M_{\rho(c_{n}),\rho(c_{n})}$ is an operator of rank $n$. It also follows from Lemma \ref{rep}  that   $M_{\rho(c),\rho(c)}-M_{\rho(c_{n}),\rho(c_{n})}=M_{\rho(c-c_{n}),\rho(c-c_{n})}$. Hence, $$a_{n}(M_{c, c})=a_{n}(M_{\rho(c),\rho(c)})\leq \|\rho(c-c_{n-1})\|^{2}\leq \lambda_{n}^{2},$$ 
and so $\rma(\Phi)\in \ell_p$.

 Assume that there exist a faithful representation $\pi : \cl
A\rightarrow \mathbf{B}(\cl H)$ and  elements $a_i,b_i\in
\cl {K}(\cl A)$ for $i=1, \dots, n$,  such that $s(\pi(a_i)), \; s(\pi(b_i))\in\ell_p$ , $i=1, \dots, n$,
and $\Phi=\sum_{i=1}^nM_{a_i, b_i}$. 
We have $\Phi(p_k)=cp_kc=\sum_{i=1}^n a_ip_k b_i$.
Hence $\pi(c)\pi(p_k)\pi(c)=\sum_{i=1}^n \pi(a_i) \pi(p_k) \pi(b_i)$
and by continuity 
\begin{equation}\label{ppp}
\pi(c)P \pi(c)=\sum_{i=1}^n \pi(a_i) P \pi(b_i).
\end{equation}
where $P=\sum_{j=1}^\infty \pi(r_j)$  is the sot-limit of the sequence  $ (\pi(p_{k}))_{k=1}^{\infty}$.

It follows from  (\ref{ppp}) that $\pi(c)P \pi(c) \in \cl S_{p/2}$. 
 On the other hand 
 $$\pi(c)P \pi(c) =\sum_{j=1}^\infty\lambda_j^2 \pi(r_{j}).$$
It follows that 
$(\lambda_j^2)\in \ell_{p/2}$
and so $(\lambda_j)\in \ell_{p},$
a contradiction.

\end{proof}

 We note the following corollary of 
Theorem \ref{th_main}.

\begin{corollary}
Let $\cl A$ be a $C^*$-algebra such that $\cl K(\cl A)$
 has finite spectrum, $\mathfrak{i}$ be a stable Calkin space and $\mathrm{s}$
be an additive $s$-function. Let $\Phi$ be an elementary operator on $\cl A$ such that $\mathrm{s}(\Phi) \in \fri$. Then $\Phi$ is a linear combination of four positive elementary operators 
$\Phi_{j}$,  $j=1,2,3,4$ such that   $\mathrm{s}(\Phi_{j}) \in \fri$ for every $j=1,2,3,4.$

\end{corollary}
\begin{proof}
 By assertion (2) of Theorem \ref{th_main},  there exist a
representation $\sum_{i=1}^{m} M_{a_i,b_i}$, 
$a_i, b_i\in \cl A,\;i=1,\dots,m$, of $\Phi$ and a faithful
representation $\pi$ of $\mathcal{A}$ such that $s(\pi(a_i)),
s(\pi(b_i)) \in \mathfrak{i}$, $i=1,\dots,m$. Let
$\Phi^{\pm}(x) = \frac{1}{4}\sum_{i=1}^m (a_i\pm b_i^*)x(a_i^*\pm b_i)$ and
$\Psi^{\pm}(x) = \frac{1}{4}\sum_{i=1}^m (a_i\pm ib_i^*)x(a_i^*\mp ib_i).$
Clearly, all operators $\Phi^{\pm}$, $\Psi^{\pm}$ are positive. By assertion (1) of 
Theorem \ref{th_main},  $\mathrm{s}(\Phi^{\pm}), \mathrm{s}(\Psi^{\pm})\in \fri$. A straightforward
verification shows that $\Phi = \Phi^+ - \Phi^- + i(\Psi^+ -
\Psi^-)$. The proof is complete.
\end{proof}
We close this section by proving  a result which may be viewed as a quantitative  
version of a result of Ylinen \cite{y}.

\begin{theorem}\label{th_main_ver}
Let $\mathcal{A}$ be a $C^{*}$-algebra, $a \in \cl A$ and  $\mathfrak{i}$  be a Calkin space. Assume that 
$\mathrm{d}(M_{a,a}) \in \mathfrak{i}$. Then   $s(\rho(a))^2\in \mathfrak{i}$ where $(\rho,\cl H)$ is the reduced atomic
representation of $\ca$
\end{theorem}
\begin{proof}
Since $\mathrm{d}(M_{a,a}) \in \mathfrak{i}$ the operator $M_{a,a}$ is compact and it follows from \cite{y} that $\rho(a)$ is compact.

Let $(\rho,\cl H)=\left(\bigoplus_{i \in I} \rho_i,\;\;
\bigoplus_{i \in I}\cl {H}_i\right)$. Set  $\cl C =
\oplus_{i\in I}\mathbf  B(\cl H_{i})$.

Let $\Phi: \cl C\rightarrow \cl C$ be the operator  defined by $\Phi(X)=\rho(a)X\rho(a)$.
 Since $\overline{\cl \rho(\cl A)}^{\mathrm{wot}} =
\cl C$, Lemma \ref{l_kap} implies that $\mathrm{d}(\Phi)\leq
\mathrm{d}(M_{\rho(a),\rho(a)})$ and so  $\mathrm{d}(\Phi) \in
\mathfrak{i}$. 
 
Let $\rho(a)=UA$ be the polar decomposition of $\rho(a)$ and $A= \sum_{k=1}^{\infty}\lambda_k
e_k^*\otimes e_k$ be a Schmidt expansion of $A$. Define $\alpha: \ell_{\infty}\rightarrow \cl C$  by $\alpha((x_{l})_{l=1}^\infty)=\sum_{l=1}^{\infty} x_l e_l^*\otimes e_l$  and $\beta: \cl C \rightarrow \ell_{\infty}$  by $\beta(X)= (\langle Xe_{l},e_{l}\rangle)_{l=1}^\infty $.    Consider the map $\Psi: \ell_{\infty}\rightarrow \ell_{\infty}$ defined by $\Psi((x_{l})_{l=1}^\infty)=\beta(U^{*}\Phi(\alpha((x_{l})_{l=1}^\infty)U^{*}))$.
Since $\alpha$  and $\beta$  are contractions we have  $\mathrm{d}(\Psi)\leq \mathrm{d}(\Phi)$ and so $\mathrm{d}(\Psi)\in\mathfrak{i}$.
A direct calculation
shows that
$\Psi((x_l)_{l=1}^\infty)= (\lambda_{l}^2 x_{l})_{l=1}^\infty.$
It follows
 \cite[Theorem 11.11.3]{pi} that
$\mathrm{d}(\Psi)=(\lambda_l^2)_{l=1}^\infty$. Hence,
$s(A)^2\in\mathfrak{i}$.
\end{proof}

\end{document}